\documentclass[12pt]{article}

\usepackage{graphicx}    
\usepackage{amsmath} 
\usepackage{amsfonts} 
\usepackage{authblk}
\usepackage{xcolor}
\usepackage{amssymb, amsthm} 
\usepackage{booktabs}
\usepackage{hyperref}    
\usepackage{geometry}    
\usepackage{caption}     
\usepackage{cite}
\usepackage{algorithm}
\usepackage{algorithmicx}
\usepackage{algpseudocode}
\usepackage{pgfplots}
\usepackage{bm}
\usepackage{float}
\usetikzlibrary{patterns,positioning}
\usepackage{tikz}
\usetikzlibrary{decorations.markings}
\usetikzlibrary{matrix}
\usepackage{pgfplots}
\usetikzlibrary{shapes, arrows.meta, patterns}
\usepackage{mathrsfs}

\usepackage{subcaption}

\newcommand{\bfa}[1]{\boldsymbol{#1}} 			
			
\newcommand{\bfeps}{\boldsymbol{\epsilon}}

\newcommand{\Div}{\text{Div}}     				%

\DeclareMathAlphabet{\mathpzc}{OT1}{pzc}{m}{it}

\DeclareMathOperator*{\argmin}{arg\, min}

\newcommand{\bfx}{\boldsymbol{x}}	
\newcommand{\bfX}{\boldsymbol{X}}	
	
\newcommand{\bfT}{\boldsymbol{T}}		
\newcommand{\bfzero}{\boldsymbol{0}}

\newcommand{\Th}{\mathscr{T}_{h}}

\newcommand{\Fh}{\mathscr{F}_{h}}

\newcommand{\refface}{\hat{\tau}} 
\newcommand{\refcoord}{\hat{\mathbf{z}}}
\newcommand{\nodes}{\mathcal{I}_h}
\newcommand{\nodecoord}{\mathbf{z}}
\newcommand{\basis}{\zeta}

\newcommand{\nodesDboundary}{\mathcal{I}_{D,h}}
\newcommand{\FESpace}{\mathcal{S}_h}
\newcommand{\FESpaceHomD}{\mathcal{S}_{d,h}}
\newcommand{\FESpaceInhomD}{\mathcal{S}^n_{g,h}}
\newcommand{\FESpaceCrack}{\mathcal{S}^n_{cr,h}}
\newcommand{\FESpaceCrackSimple}{\mathcal{S}_{cr,h}}
\newcommand{\FESpaceInhomDSimple}{\mathcal{S}_{g,h}}
\newcommand{\FacesCR}{\mathscr{F}_h^{CR}}
\newcommand{\CrackRegion}{C_h}
\newcommand{\tolCR}{\Xi_{CR}}
\newcommand{\phasefield}{\psi_h}

\usepackage{etoolbox}

\newtheorem*{dwf}{Discrete formulation}

\newtheorem{remark}{Remark}

\usepackage{etoolbox} 
\makeatletter
\patchcmd{\blfootnote}{\@mpfootnotetext}{\@gobble}{}{}
\makeatother

\newtheoremstyle{mypropositionstyle} 
  {\topsep}                     
  {\topsep}                     
  {\itshape}                    
  {}                            
  {\bfseries}                   
  {.}                           
  {.5em}                        
  {\thmname{#1}\thmnumber{ #2}\thmnote{ (#3)}}

\theoremstyle{conjecturestyle}

\newtheorem{proposition}{Proposition}

\newenvironment{form1}[1]{\par\medskip\noindent\textbf{Formulation 1}\label{#1}\itshape}{\par\medskip}

\title{An \textsf{AT1} phase-field framework for quasi-static anti-plane shear fracture: Unifying $\xi$-based adaptivity and nonlinear strain energy density function\footnote{This work is dedicated to the cherished memory of Professor K. R. Rajagopal (Professor of Mechanical Engineering, Texas A\&M University, College Station, Texas-USA). His seminal contributions and invaluable mentorship were not only instrumental in shaping this research but also profoundly influenced the author SMM.}}
\author[1]{Maria P. Fernando}
\author[1,*]{S. M. Mallikarjunaiah}

\affil[1]{Department of Mathematics \& Statistics, Texas A\&M University-Corpus Christi, TX- 78412, USA}
\affil[*]{Corresponding author}
\affil[ ]{\textit{E-mail addresses:} \texttt{fpieo@islander.tamucc.edu} (M.P. Fernando), \texttt{M.Muddamallappa@tamucc.edu} (S.M. Mallikarjunaiah)}
\date{}

\begin{document}

\maketitle  
\begin{abstract}
This work introduces a novel \textsf{AT1} phase-field framework for simulating quasi-static anti-plane shear fracture in geometrically linear elastic bodies. A key feature of this framework is the unification of $\xi$-based local mesh adaptivity---where $\xi$ represents the characteristic length of the damage zone---and an algebraically nonlinear strain energy density function. A modified Francfort-Marigo energy functional, together with its Ambrosio-Tortorelli-type regularization, is hereby proposed to address challenges within the framework of nonlinearly constituted materials.  Such a nonlinear constitutive model is crucial for circumventing the physically inconsistent crack-tip strain singularities that arise in classical linear elastic fracture mechanics. We dynamically optimize $\xi$ throughout the simulation, significantly enhancing the computational efficiency and accuracy of numerically approximating the local minimizers of the Ambrosio-Tortorelli (\textsf{AT1})-type phase-field model. The proposed regularization for the total energy functional comprises three distinct components: a nonlinear strain energy, an evolving surface energy, and a linear-type regularization term dependent on the length scale of the damage zone. Variational principles applied to this novel energy functional yield a coupled system of governing second-order quasilinear partial differential equations for the mechanics and phase-field variables. These equations are subsequently discretized using the conforming bilinear finite element method. The formulation is underpinned by four crucial parameters: two are integral to the nonlinear strain energy function, while the other two serve as penalty parameters. These penalty parameters are asymptotically calibrated and rigorously utilized in the numerical simulations. Our results demonstrate that this spatially adaptive approach leads to enhanced mesh adaptivity, ensuring the robust convergence of the numerical solution. Numerical examples show that the proposed adaptive model significantly outperforms standard phase-field methods by providing a more accurate representation of fracture propagation while substantially reducing computational costs. Furthermore, the employment of this adaptive strategy allows for the use of a vastly larger regularization length parameter throughout the entire computation compared to non-adaptive counterparts.
\end{abstract}

\vspace{.1in}

\textbf{Key words.} Quasi-static crack; anti-plane shear loading; spatially adaptive; Ambrosio-Tortorelli energy functional; Continuous Galerkin finite element method; Damage zone length

\section{Introduction}
The phenomenon of cracks and fractures in materials, particularly in elastic and porous media, stands as a cornerstone of contemporary research in applied mathematics and engineering \cite{gou2015modeling,ghosh2025finite,manohar2025convergence,rajagopal2011modeling,Mallikarjunaiah2015,ferguson2015,yoon2024finite,yoon2022finite,yoon2022CNSNS,silling2000reformulation,vasilyeva2024generalized,manohar2024hp}. This field, broadly known as fracture mechanics, is paramount for guaranteeing the safety and longevity of engineered systems, spanning from consumer goods to critical infrastructure. A deep understanding of crack initiation and propagation is indispensable for mitigating catastrophic failures in structures such as bridges, aerospace components, and pipelines. Fundamentally, this knowledge facilitates the design of more robust materials and components, enabling precise predictions of structural service life, thereby averting potential human casualties and substantial economic losses.

Modeling crack behavior is inherently challenging because it's so complex to describe the entire process with a fixed set of equations. Several studies have tackled crack evolution using various methods. Some explicitly define crack-tip movement, like the extended finite element method \cite{fries2010extended} and the generalized finite element method \cite{babuvska2012stable}, among others. More recently, variational regularization approaches, often called phase-field approximation methods, have gained prominence for studying both quasi-static \cite{francfort1998revisiting,bourdin2008variational,bourdin2000numerical,yoon2021quasi,lee2022finite,burke2013adaptive,burke2010adaptive} and dynamic \cite{larsen2010models,larsen2010existence,bourdin2011time,dal2011existence,dal2016existence,borden2012phase} crack evolution. phase-field fracture models eliminate the need to explicitly track sharp crack discontinuities. Instead, they use a regularized representation where a continuous, scalar phase-field variable, typically denoted by $v$, describes the material's integrity. This variable smoothly transitions from a value representing intact material to a value representing the fully fractured state, effectively creating a diffuse interface or transition zone of finite width. A primary strength of this methodology lies in its intrinsic ability to predict complex fracture phenomena without needing ad-hoc criteria. The governing equations naturally capture a crack's entire lifecycle---from its initial nucleation and subsequent propagation to complex geometric evolutions like branching, kinking, and the formation of arbitrary curvilinear pathways. Consequently, this approach sidesteps the cumbersome procedures inherent to traditional fracture mechanics. It eliminates the need for complex post-processing of crack-tip fields to compute criteria like stress intensity factors. Furthermore, it avoids the computationally expensive and algorithmically challenging task of repeatedly remeshing the domain to conform to the evolving crack geometry.

The foundational theory of brittle fracture, pioneered by Francfort and Marigo \cite{francfort1998revisiting}, ingeniously reframed the intricate problem of crack propagation within a variational energy minimization framework. At its core, this model postulates that cracks initiate and evolve to continuously minimize the system's total energy. This formulation is rooted in an energetic principle that balances two competing physical processes: the driving force of fracture, manifested as the release of stored elastic strain energy from the bulk material, and the resistance offered by fracture energy, representing the energy consumed to create new crack surfaces. Direct computational implementation of this principle is challenging due to the mathematical and algorithmic complexities associated with tracking an evolving, geometrically sharp discontinuity. To overcome this hurdle and enhance computational tractability, the sharp-crack topology is typically regularized using the \textit{Ambrosio-Tortorelli} functional. This well-established phase-field approximation replaces the distinct crack surface with a diffuse damage zone, effectively smoothing the problem. This powerful mathematical technique approximates the sharp discontinuity by introducing a continuous "phase-field" variable, which smoothly transitions between the undamaged and fully broken states of the material. Even with this regularization, the direct, monolithic solution of the resulting full coupled, nonlinear Euler-Lagrange equations presents significant stability and convergence problems. Therefore, these systems are typically addressed using sophisticated numerical approaches, such as staggered solution schemes, line-search algorithms, or other advanced methods to ensure robust and accurate results \cite{gerasimov2016line,miehe2010phase,amor2009regularized,natarajan2019phase}.

While the length scale parameter is a critical component for regularizing sharp crack topologies in phase-field models (PFMs) of brittle fracture, its introduction presents several significant challenges. A primary issue concerns its physical interpretation; it remains debated whether it is a purely numerical regularization parameter or a genuine material property representing the material's microstructure \cite{wu2018length,bourdin2008variational}. In many classical PFMs, the length scale directly governs the predicted tensile strength of the material, creating an undesirable dependence of a physical property on a numerical parameter \cite{miehe2010thermodynamically}. This dependency is intrinsically linked to computational cost. Accurately resolving the diffuse damage zone necessitates a mesh size significantly smaller than the length scale. Consequently, approximating a sharp crack by choosing a very small length scale leads to prohibitive computational expense, particularly for large-scale, three-dimensional simulations \cite{borden2012phase}. Furthermore, standard PFMs often suffer from spurious damage initiation at stress concentrations, where damage may evolve prematurely before the global failure criterion is met \cite{amor2009regularized}. This artifact, which also depends on the length scale, can compromise the predicted stiffness and peak load of the structure.

In response to these limitations, significant research efforts aim to develop more advanced formulations---such as those employing hybrid or unified degradation functions---that decouple the regularization length scale from the physical fracture response \cite{wu2017unified}. In recent work \cite{phansalkar2022spatially,phansalkar2023extension}, a model extends the classical Ambrosio-Tortorelli (\textsf{AT2}) framework by including the damage length scale parameter as a third unknown variable in the global energy minimization. However, the foundational \textsf{AT2} model possesses a well-documented limitation \cite{pham2011gradient}. Specifically, its formulation is nonlinear in elastic variation prior to fracture, which induces a non-physical artifact of early material degradation, causing the model to predict a loss of stiffness at low stress levels where the material should behave purely elastically.

Understanding how cracks behave in materials is crucial for predicting when and where failure might occur, as well as for assessing the forces at play locally. Analyzing the deformation near a crack provides vital insights into a material's capacity to stretch and deform before it breaks or suffers microscopic damage \cite{broberg1999}. Moreover, considering the strain energy density, which combines both stress and strain, offers a more comprehensive way to characterize the energy absorbed and released at the crack tip. This provides a robust tool for predicting crack initiation and growth, even under complex loading conditions \cite{yoon2021quasi}. For materials with directional properties, ignoring the impact of uneven loading on these crucial crack-tip parameters could severely limit our ability to predict their behavior, potentially jeopardizing the safety and integrity of components in high-performance applications.


Accurately mapping stress and strain distributions around features like notches, slits, or holes is a fundamental challenge in both engineering and theoretical mechanics. Traditionally, the analytical framework for understanding these stress concentrations has relied on linearized elasticity theory \cite{Inglis1913, lin1980singular, love2013treatise, murakami1993stress}. However, a major limitation of this classical approach is its prediction of infinite strain at the tips of such discontinuities. This is physically unrealistic and stems from its simplified approximation of finite deformations. This inconsistency has prompted extensive research into developing more realistic material models \cite{gurtin1975, sendova2010, MalliPhD2015, ferguson2015, zemlyanova2012, WaltonMalli2016, rajagopal2011modeling, gou2015modeling}, often coupled with advanced numerical methods like collocation techniques. Despite these efforts, a significant obstacle remains: balancing model accuracy with computational efficiency and ease of experimental validation \cite{broberg1999}. Many proposed theoretical improvements, while more realistic, often come with high computational costs or are difficult to verify experimentally.


Furthermore, linear elastic fracture mechanics (LEFM), despite its widespread use for modeling crack behavior, has inherent limitations. Beyond the well-known strain singularity, LEFM also predicts an unrealistic blunt crack profile and the problematic interpenetration of crack faces, particularly in materials joined together. The issue of crack-tip singularity isn't fully resolved even within various nonlinear elasticity frameworks, as shown in studies like \cite{knowles1983large} and models incorporating specific constraints, such as the bell constraint model \cite{tarantino1997}. These ongoing challenges raise a key question: Can nonlinear algebraic material models effectively alleviate the crack-tip strain singularity, even if some stress singularity persists? This question is a primary motivation for the current investigation.


Traditional elasticity theories, such as those by Cauchy and Green, often fall short when describing material behavior under extreme deformations, especially at crack tips where unrealistic singularities can appear. To overcome these limitations, Rajagopal and his collaborators developed a generalized framework for elasticity \cite{rajagopal2003implicit,rajagopal2007elasticity, rajagopal2007response,rajagopal2009class,rajagopal2011non,rajagopal2011conspectus,rajagopal2014nonlinear,rajagopal2018note}. This innovative approach, \textit{Rajagopal's theory of elasticity}, uses implicit constitutive models grounded in robust thermodynamic principles. Within this framework, an elastic body, by definition non-dissipative, is characterized by implicit relationships connecting the Cauchy stress and deformation gradient tensors \cite{bustamante2018nonlinear}. A particularly powerful aspect of Rajagopal's theory is its ability to produce a hierarchy of 'explicit' nonlinear relationships, allowing linearized strain to be expressed as a nonlinear function of stress. Crucially, a specific subclass of these implicit models has a unique 'strain-limiting' property: they can represent linearized strain as a uniformly bounded function across the entire material domain, even under significant stress conditions. This characteristic makes these models exceptionally well-suited for investigating crack and fracture behavior in brittle materials \cite{rajagopal2011modeling, gou2015modeling, Mallikarjunaiah2015, MalliPhD2015}, offering a path toward analyzing both quasi-static and dynamic crack evolution. The usefulness of these strain-limiting models has been demonstrated through numerous studies that have revisited and provided new insights into classical elasticity problems \cite{kulvait2013,rajagopal2018bodies,bulivcek2014elastic,erbay2015traveling,zhu2016nonlinear,csengul2018viscoelasticity,itou2018states,itou2017contacting,yoon2022CNSNS}. Their versatility in explaining the mechanical behavior of a wide range of materials, especially regarding crack and fracture phenomena, is a significant advantage. Recent research, for example, has shown that formulating quasi-static crack evolution problems within this strain-limiting framework can predict complex crack patterns and even increased crack-tip propagation velocities \cite{lee2022finite, yoon2021quasi}.

This study marks a significant advancement in the field by integrating the damage length scale parameter as a third primary variable into the well-established Ambrosio-Tortorelli (\textsf{AT1}) phase-field framework. Beyond this, we propose a revised Francfort-Marigo total energy function, specifically by incorporating a nonlinear strain-energy density function within the principles of strain-limiting theories of elasticity. A key contribution of this work is the development of a revised AT1-type regularization, which further refines the modeling of fracture processes. These novel inclusions collectively enhance the framework's ability to accurately capture complex material behaviors under various loading conditions. Consequently, the total energy functional in this framework is comprehensively defined by three distinct components. It incorporates a nonlinear strain energy density function, which accounts for the material's elastic response, particularly under significant deformations. Additionally, an \textsf{AT1}-type crack surface energy term captures the energetic cost associated with the formation and propagation of cracks. Finally, a crucial third term, which is a function of the damage length scale parameter, is included. This latter term provides the necessary mechanism for dynamically determining and optimizing the characteristic length associated with damage localization. The model's total energy functional, comprising three distinct terms, is simultaneously minimized with respect to all three variables. A key innovation of our methodology is the use of this dynamically optimized and spatially varying length parameter. It serves as a natural and physically motivated indicator, effectively guiding local adaptive mesh refinement. To solve the coupled system of partial differential equations, which are derived from the Euler-Lagrange equations of the total energy functional, we employ a conforming Continuous Galerkin finite element method. Our results reveal that the optimal length parameter, whether treated as a global or locally varying field, is considerably larger than values typically used in standard regularized phase-field models. This larger parameter offers a twofold advantage: it reduces computational cost by permitting a coarser mesh discretization, and it acts as a crucial factor in achieving robust numerical convergence for the approximation.

This paper is structured as follows: Section \ref{intro_model} introduces the governing model for quasi-static brittle fracture and presents a modified energy functional designed to compute an optimal regularization length parameter. We formulate two distinct variants: one determines a single, optimal uniform length parameter for the entire domain, while the other treats the length scale as a spatially varying field, solved concurrently with the displacement and phase-fields. Section \ref{fem} details the finite element discretization of this proposed three-field problem, including both its continuous variational form and final discrete algebraic formulation. Section \ref{rd} focuses on numerical implementation and results, discussing key numerical aspects, parameter calibration, and the staggered iterative algorithm used for computations. Subsequently, we present and analyze numerical results from benchmark tests involving anti-plane shear loading. Finally, Section \ref{conclusions} summarizes the key findings and outlines potential directions for future research.

\section{Mathematical model for quasi-static fracture in strain-limiting solids}\label{intro_model}
Consider a material body represented by a closed and bounded Lipschitz domain $\Lambda \subset \mathbb{R}^2$ in its reference configuration. The Lipschitz continuous boundary $\partial \Lambda$ is partitioned into two disjoint sets: $\partial \Lambda = \overline{\Gamma_{D}} \cup \overline{\Gamma_{N}}$, where $\Gamma_{D}$ is a non-empty Dirichlet boundary and $\Gamma_{N}$ is a Neumann boundary. Let $\bfx = (X_1, X_2)$ and $\bfx = (x_1, x_2)$ denote typical points in the reference and deformed configurations, respectively. The displacement field is given by $\mathbf{u} \colon \Lambda \to \mathbb{R}^2$, where $\mathbf{u} = \bfx - \bfX$.

For the remainder of this paper, we adopt standard notations from Lebesgue and Sobolev spaces, as detailed in references like \cite{ciarlet2002finite,evans1998partial}. The space of Lebesgue integrable functions $L^{p}(\Lambda)$ for $p \in [1, \infty)$ is defined, with $L^{2}(\Lambda)$ specifically denoting the space of square-integrable functions equipped with the inner product $\left( v, w \right) := \int_{\Lambda} v \, w \, d\bfx$ and induced norm $\| v \|_{L^2(\Lambda)} := \left(v, v \right)^{1/2}$. $C^{m}(\Lambda)$ represents the space of $m$-times continuously differentiable functions on $\Lambda$, for $m \in \mathbb{N}_0$. The standard Sobolev space $H^{1}(\Lambda)$ is defined as:
$$
H^{1}(\Lambda) := W^{1, 2}(\Lambda) = \left\{ v \in L^{2}(\Lambda) \; : \; \partial_j v \in L^{2}(\Lambda) \text{ for } j \in \{ 1, 2 \} \right\},
$$
and is endowed with its standard inner product and induced norm. The closure of $C_0^{\infty}(\Lambda)$ (the space of infinitely differentiable functions with compact support) under the $H^1(\Lambda)$ norm is denoted $H_0^{1}(\Lambda)$, i.e., $H_0^{1}(\Lambda) = \overline{C_0^{\infty}(\Lambda)}^{\| \cdot \|_{H^1}}$. A function $f \in L^{1}(\Lambda)$ is said to have bounded variation if \cite{BuOrSue10, BuOrSue13}:
$$
\sup \left\{ \int_{\Lambda} f \, \operatorname{div} \boldsymbol{\psi} \, d\bfx \; : \; \boldsymbol{\psi} \in \mbox{C}_{0}^{1}\left( \Lambda; \, \mathbb{R}^{N} \right), \; |\boldsymbol{\psi}| \leq 1 \right\} < \infty.
$$
The space of such functions is denoted by $\operatorname{BV}(\Lambda)$. Functions in $\operatorname{BV}(\Lambda)$ may exhibit discontinuities, which are reflected in their distributional gradient. The space of special functions of bounded variation, $\operatorname{SBV}(\Lambda)$, introduced in \cite{AmTo92}, is a subset of $\operatorname{BV}(\Lambda)$. For any $f \in \operatorname{SBV}(\Lambda)$, its distributional derivative is given by:
$$
Df = \nabla f \mathcal{L}^{N} + \left( f^{+}( \bfx) - f^{-}(\bfx) \right) \otimes \boldsymbol{\nu}_{f}( \bfx) \, \mathcal{H}^{N-1} \lfloor S(f),
$$
where $\nabla f$ is the approximate gradient of $f$, $S(f)$ is the jump-set of $f$, $\boldsymbol{\nu}_f$ is the unit normal to the jump-set, $f^+$ and $f^-$ are the approximate trace values of $f$ on either side of $S(f)$, and $\mathcal{L}^{N}$ and $\mathcal{H}^{N-1}$ are the $N$-dimensional Lebesgue measure and $(N-1)$-Hausdorff measure, respectively.

\subsection{A new class of response relations}

This work investigates a subclass of general implicit constitutive relations for elastic bodies, as discussed in \cite{rajagopal2003implicit,rajagopal2007elasticity,rajagopal2011non,rajagopal2014nonlinear,rajagopal2007response}. Classical Cauchy elasticity posits a constitutive relationship where the Cauchy stress $\bfT$ is an isotropic tensor-valued function of the deformation gradient $\mathbf{F}$:
\begin{equation}
\mathbf{0} = \mathcal{F}(\mathbf{B}, \bfT).   
\end{equation}
where $\mathcal{F}$ is assumed to be an isotropic function \cite{gou2015modeling}. 
A specific subclass introduced by Rajagopal \cite{rajagopal2007elasticity} is:
\begin{equation}\label{SL1-rephrased}
\mathbf{B} := \mathcal{F}( \bfT).
\end{equation}
If, for this class of models, there exists a constant $M > 0$ such that 
$$
\sup_{\bfT \in \operatorname{Sym}} |\mathcal{F}( \bfT)| \leq M,
$$ 
these relations are termed \textit{strain-limiting constitutive theories} \cite{MalliPhD2015,Mallikarjunaiah2015,gou2015modeling,Mallikarjunaiah2015,MalliPhD2015,itou2017nonlinear,yoon2022finite,manohar2024hp}.  Standard linearization procedures from classical elasticity lead to:
\begin{equation}\label{sl_model}
\boldsymbol{\varepsilon} = \beta_1 \mathbf{I} + \beta_2 \bfT + \beta_3  \bfT^2,
\end{equation}
 Finally, the system of partial differential equations governing the problem within the framework of algebraically nonlinear strain-limiting  elasticity is:
\begin{subequations}
\begin{align}
-\nabla \cdot \bfT &=\mathbf{0}, \quad \text{and} \quad \bfT= \bfT^T, \label{equilib:eq-rephrased} \\
\boldsymbol{\varepsilon} &= \Psi_{0}\left( \operatorname{tr} \, \bfT, \; \| \bfT \| \right) \mathbf{I} + \Psi_{1}\left( \| \bfT \| \right) \bfT, \quad \Psi_{0}\left( 0, \cdot \right)=0, \label{eqn:main-rephrased} \\
\operatorname{curl} \, \operatorname{curl} \, \boldsymbol{\varepsilon} &=\mathbf{0}, \quad \text{and} \quad
\boldsymbol{\varepsilon} = \frac{1}{2} \left( \nabla \mathbf{u} + \nabla \mathbf{u}^{T} \right).
\end{align}
\end{subequations}
In these relationships, $\Psi_{0} \colon \mathbb{R} \times \mathbb{R}_{+} \to \mathbb{R}$ and $\Psi_{1} \colon \mathbb{R}_{+} \to \mathbb{R}$ are scalar functions of the invariants of the Cauchy stress tensor.  The Equation~\ref{eqn:main-rephrased} represents a \textit{specialized subclass} within a broader category of material models, specifically those conforming to the general class described by Equation~\ref{sl_model}. This particular subclass offers an attractive framework for understanding material behavior under certain conditions. Over the years, a diverse range of constitutive relations has been developed to describe the mechanical response of various materials accurately. For elastic solids, numerous models have been proposed and extensively studied, as evidenced in \cite{itou2018states,itou2017nonlinear,kulvait2019state,bustamante2021new,bustamante2020novel}. These contributions have significantly advanced our understanding of how elastic materials deform and transmit stress. Furthermore, the complexities of  \textit{porous elastic solids} have also been a significant area of research, with constitutive models introduced in \cite{gou2025computational,murru2021stress,pruuvsa2022pure,rajagopal2021implicit,yoon2024finite,itou2025nonlinear}. These models are crucial for analyzing materials with void spaces, which exhibit unique mechanical properties due to the interaction between the solid matrix and the pores. The continuous development and refinement of such constitutive laws are essential for accurate material characterization and for predicting their behavior in engineering applications.
 
\subsection{Anti-plane shear problem (Mode-III)}

This study investigates the quasi-static evolution of a crack under anti-plane shear (Mode-III) loading. This specific loading condition is characterized by a planar deformation, meaning that all kinematic quantities---the displacement vector $\bfa{u}(\bfa{x}, t)$, stress tensor $\bfT(\bfa{x},t)$, and strain tensor $\bfeps(\bfa{x},t)$---are exclusively dependent on the in-plane coordinates $x_1$ and $x_2$. Consequently, displacements within the body are confined to a single direction, specifically the $x_3$-direction, which itself depends solely on $x_1$ and $x_2$. This results in a displacement field where the only non-zero component is:
\begin{equation}\label{eq:disp_vector}
\bfa{u}(x_1,\, x_2,\, t) = \left( 0, 0, u( x_1,\, x_2,\, t) \right).
\end{equation}
For a linear, isotropic, and homogeneous material, the foundational constitutive relationship is given by Hooke's Law:
\begin{equation}
\bfT = 2 \, \mu \, \bfeps,
\end{equation}
where $\mu$ represents the shear modulus of the material. However, within the framework of the strain-limiting theory of elasticity, this relationship is generalized to:
\begin{equation}
\bfeps = \Psi_{1}\left( \| \bfT \| \right) \bfT. \label{eqn:main2}
\end{equation}
To formulate the governing equation, we define the stress components through the derivatives of a scalar-valued Airy's stress function, $\Phi$. This function inherently satisfies the balance of linear momentum ($\Div \, \bfT = \bfzero$). By combining the strain compatibility condition with the nonlinear constitutive relationship presented in Equation \eqref{eqn:main2}, we derive a second-order quasi-linear partial differential equation:
\begin{equation}\label{pde:nlin1}
- \nabla \cdot \left( \Psi_{1} \left( \| \nabla \Phi \| \right) \; \nabla \Phi \right) =0,
\end{equation}
where $\| \cdot \|$ is the norm of an vector field.  For the remainder of this paper, we adopt a specific form for the constitutive function $\Psi_{1}$:
\begin{equation}\label{eq:Psi}
\Psi_{1} (\|\bfT\|) = \frac{1}{2 \, \mu \left( 1 + \beta^\alpha \, \| \bfT \|^{2\alpha} \right)^{1/\alpha}}.
\end{equation}
In this formulation, $\beta \geq 0$ and $\alpha >0$ are the modeling parameters, allowing for flexibility in describing material behavior.

\begin{remark}
When $\beta$ approaches zero ($\beta \to 0$), the model described by Equation \eqref{eq:Psi} recovers the linear elastic constitutive relationship, signifying elastic behavior under small deformations. Conversely, a large value of $\beta$ ($\beta \gg 1$) leads to a nonlinear constitutive relationship, even within the infinitesimal strain regime, indicating significant deviation from linearity. Furthermore, the function $\Psi_1$ possesses properties of invertibility and monotonicity, which are crucial for generating both hyperelastic and, in certain cases, non-hyperelastic relations, as supported by previous research \cite{mai2015strong,mai2015monotonicity,rajagopal2007response,rajagopal2011conspectus}.
\end{remark}
\begin{remark}
A key characteristic of the chosen monotone decreasing function $\Psi_{1} (r)$ in Equation \eqref{eq:Psi} is its limiting behavior:
\begin{equation}\label{EqLimit}
\lim_{r \to \infty} r \Psi_{1} (r) \to \frac{1}{2 \mu \beta}.
\end{equation}
An important consequence of this limit, as shown in Equation \eqref{EqLimit}, is that as the magnitude of stress approaches infinity ($\| \bfT \| \to \infty$), the resulting strains remain uniformly bounded throughout the material body. This holds true even in regions of high stress concentration, such as near crack tips or re-entrant corners, which is a crucial aspect for both numerical stability and physical realism.
\end{remark}
In view of Equation~\eqref{eq:Psi}, the nonlinear PDE (Equation~\eqref{pde:nlin1}) now takes the form:
\begin{equation} \label{pde:mech}
- \nabla \cdot \left( \frac{\nabla \Phi}{2 \, \mu \left( 1 + \beta^\alpha \; \|\nabla \Phi \|^{2\alpha} \; \right)^{1/\alpha}} \right) = 0.
\end{equation}

\begin{remark}
This quasi-linear partial differential equation (Equation \eqref{pde:mech}), particularly with specific selections for the modeling parameters $\beta$ and $\alpha$, exhibits notable similarities to the classical minimal surface optimality problem in variational mechanics. Previous work, particularly in \cite{bulivcek2015existence,bulivcek2014elastic,bulivcek2015analysis}, has demonstrated the existence of a continuous solution for the anti-plane shear problem within a non-convex v-notch domain for $\beta \in (0, \, 2)$ and $\alpha > 0$, utilizing dual variational methods.
\end{remark}

\begin{remark}
Despite exhibiting significant nonlinear behavior at low strains (e.g., within 1\%), many materials—such as gum metal \cite{kulvait2019state}, titanium alloys commonly used in orthopedic applications \cite{tian2015nonlinear,devendiran2017thermodynamically}, and Ti-30Nb-10Ta-5Zr (TNTZ-30) alloys \cite{hao2005super,saito2003multifunctional}—have traditionally been analyzed using linear elastic models. This conventional approach, which assumes a linear relationship between infinitesimal strain and Cauchy stress, is inherently problematic given the observed nonlinearity. A more fitting framework for describing such behavior is the implicit theory of elasticity \cite{rajagopal2014nonlinear}. Additionally, strain-limiting nonlinear relationships effectively characterize materials like rocks \cite{bustamante2020novel} and various rubber-like compounds \cite{bustamante2021new}. For example, a constitutive model for rubber introduced by Bustamante et al. \cite{bustamante2021new}, which employs principal stresses as its fundamental variables, has shown superior alignment with experimental data compared to the widely used Ogden's model \cite{ogden1972large}. 
\end{remark}

The problem outlined in Equation \ref{pde:mech}, concerning a non-convex domain featuring a crack or slit subjected to classical zero-traction boundary conditions on the crack faces, was investigated in \cite{manohar2024hp}. It was demonstrated in that study that a unique solution to the discrete problem, which was discretized using a local discontinuous Galerkin (DG) finite element method, is ensured for specific ranges of the modeling parameters $\beta$ and $\alpha$. Furthermore, it was observed that while the $hp$-version of the DG method yields an optimal convergence order with respect to the mesh size $h$, a suboptimal convergence rate was exhibited with respect to the polynomial degree $p$. 

The intricate problem of modeling the mechanical response of brittle materials, specifically those containing cracks or wedges within both isotropic and orthotropic solid frameworks, has been extensively investigated in prior works, including \cite{kulvait2019state,kulvait2013,yoon2022finite,ghosh2025finite,ghosh2025b}. Across these studies, it has been consistently observed that the parameters $\alpha$ and $\beta$ exert a profound and critical influence on the localized stress and strain concentrations near the crack tip. More specifically, the parameter $\beta$ is primarily identified as contributing to a toughening or crack-mitigating effect, by effectively governing the degree of strain concentration. Conversely, the parameter $\alpha$ has been found to exert a detrimental influence on the crack-tip fields, thereby directly impacting the material's susceptibility to fracture initiation and ultimate failure. It has been repeatedly demonstrated that a progressive increase in the value of $\alpha$ invariably leads to a marked intensification of both stress and strain concentrations at the crack tip.
 
\subsection{Modified Francfort-Marigo model and its regularization for quasi-static evolution}
Within the scope of this work, the domain of interest, denoted as $\Lambda:=\Lambda(t) \in \mathbb{R}^2$, is established as a smooth, open, connected, and bounded region, possessing a defined boundary $\partial \Lambda$. The presence of a crack set, $\Gamma(t)$, is also considered. It is understood that the time variable, $t \in [0, T]$, where $T$ represents the final time, influences the system exclusively through time-dependent loading conditions. A key assumption made is that the discontinuity set $\Gamma(t)$ is entirely contained within $\Lambda(t)$ and is characterized as a Hausdorff measurable set.

The current investigation is initiated by considering the Francfort-Marigo energy functional \cite{francfort1998revisiting,bourdin2008variational}, which was originally proposed within the context of linearized elasticity. This functional, expressed as:
\begin{equation}\label{eq:tenergy}
{E}(\Phi, \; \Gamma):=\int_{\Lambda \setminus \Gamma} \widehat{\mathcal{W}}(\Phi) d\bfx + \mathcal{G}_c \mathcal{H}^{d-1}(\Gamma),
\end{equation}
is understood to relate the rate of decrease in the material's strain energy to the ``cost'' associated with the creation of new crack surfaces. This formulation is entirely consistent with Griffith's criterion for fracture \cite{griffith1921phenomena}. In this expression, $\widehat{\mathcal{W}}(\cdot) \colon H^{1}(\Lambda) \to \mathbb{R}$ represents the elastic energy (which can be derived from either linear or nonlinear strain-limiting elasticity or strain limiting theories elasticity \cite{rajagopal2007elasticity,rajagopal2011modeling}). Furthermore, $\mathcal{H}^{d-1}$ denotes the Hausdorff measure, and $\mathcal{G}_c$ represents the critical energy release rate (also known as the fracture toughness) of the material.

The unilateral minimization of this total energy functional is known to yield a new equilibrium state and an updated crack set. The time-dependent minimization of this energy functional has been extensively studied for the existence of solutions, primarily through the application of variational calculus principles.  For any $\Phi \in \mbox{SBV}(\Lambda)$ and an associated discontinuity set $\Gamma$, the bulk, surface, and total energies are respectively defined as:
\begin{align}
E_{B} (\Phi) &:=\int_{\Lambda\backslash J(\Phi)}  \widehat{\mathcal{W}}(\Phi) \; d\bfx, \quad \mbox{with} \quad  \widehat{\mathcal{W}}(\Phi) = \dfrac{ \| \nabla \Phi\|^2}{ \left(1 + \beta^{\alpha} \| \nabla \Phi \|^{2\alpha}\right)^{1/\alpha}} \label{eq:bulk_energy} \\
E_S(\Gamma) &:= \mathcal{G}_c \, \mathcal{H}^{N-1} (\Gamma) \label{eq:surface_energy} \\
E(\Phi, \, \Gamma) &:= \begin{cases} 
 E_{B} (\Phi) + E_S(\Gamma) & \mbox{if} \quad  \mathcal{H}^{N-1} \left( J(\Phi) \backslash \Gamma \right)  = 0\\
+ \infty & \mbox{otherwise }  
\end{cases} \label{eq:total_energy}
\end{align}
Here, $J(\Phi)$ denotes the \textit{jump set} of the function $\Phi$.

In this paper, a central focus is placed on the study of quasi-static evolution within brittle materials when subjected to a time-varying load $g(t) \in L^{\infty}(0, \, T; W^{1, \, \infty}(\Lambda)) \cap W^{1, \, 1}(0, \, T; H^{1}(\Lambda))$. This load is assumed to be applied on a specific open subset $\Lambda_D \subset \Lambda$. For the sake of notational simplicity concerning boundary conditions, the following set is defined:
\begin{equation} \label{eq:boundary_set}
\mathcal{D}(g(t)) :=\left\{ \Phi \in \mbox{SBV}(\Lambda) \colon \Phi\arrowvert_{\Lambda_D} = g(t) \right\}.
\end{equation}
The quasi-static formulation is typically structured as a successive global minimization problem \cite{bourdin2008variational,bourdin2000numerical} to determine the equilibrium solution, as shown:
\begin{subequations}\label{eq:minimization}
\begin{align}
\Phi(t_k) &:= \argmin_{{v \in \mathcal{D}(g(t_k))}} E_B(v) + E_S(J(v) + \Gamma(t_{k-1})), \label{eq:minimization_phi} \\
\Gamma(t_k) &:= J(\Phi(t_k)) \cup \Gamma(t_{k-1}), \label{eq:minimization_gamma}
\end{align}
\end{subequations}
 
 \subsubsection{Algebraically nonlinear strain energy and Ambrosio-Tortorelli-type regularization}
A direct numerical implementation of the global minimization problem presented in Equation \eqref{eq:minimization} is confronted with several inherent challenges. These primarily stem from the irregularity of the total energy functional and the complexities associated with discretizing the unknown crack set, which inherently possesses an unknown path. A successful approach, widely adopted in the literature, involves the regularization of the total energy functional $E(\Phi, \, \Gamma)$ to enable the tracking of the crack set, thereby making it amenable to discretization methods. For given positive parameters $\kappa >0$ and $\xi >0$, the regularized elastic energy $\widetilde{\mathcal{E}} \colon H^1 \times H^1 \to \mathbb{R} \cup \{ + \infty\}$ and the regularized crack-surface energy $\mathcal{H} \colon H^1 \to \mathbb{R}$ are defined, respectively, as:
\begin{subequations}\label{eq:regularization}
\begin{align}
\widetilde{\mathcal{E}}(\Phi, \, \varphi) &:= \frac{\mu}{2} \int_{{\Lambda}} \left( (1- \kappa) \varphi^2 + \kappa  \right) \;  \mathcal{W}(\Phi, \, \varphi)  \; d \bfx,  \\
 \mathcal{H}(\varphi)&:= \dfrac{1}{c_w}\int_{{\Lambda}} \left[ \frac{(1-\varphi)}{\xi} + \xi \; \|\nabla \varphi\|^2  \right]  \; d\bfx,
\end{align}
\end{subequations}
In the above equations, the regularized strain-energy term 
\begin{equation}
\mathcal{W}(\Phi, \, \varphi) = \dfrac{ \| \nabla \Phi\|^2}{ \left(1 + \beta^{\alpha} \left( (1- \kappa) \varphi^2 + \kappa  \right)^\alpha \| \nabla \Phi \|^{2\alpha}\right)^{1/\alpha}}
\end{equation}

Consequently, the total energy of the mechanical system, as established in various works \cite{AmTo92}, is given by:
\begin{align}
E_{\xi}(\Phi, \, \varphi) &:=   \widetilde{\mathcal{E}}(\Phi, \, \varphi) + \mathcal{G}_c \;  \mathcal{H}(\varphi) \notag \\
&= \frac{\mu}{2} \int_{{\Lambda}} \left( (1- \kappa) \varphi^2 + \kappa  \right) \; \mathcal{W}(\Phi, \, \varphi)  \; d \bfx + \dfrac{\mathcal{G}_c}{c_w}\int_{{\Lambda}} \left[ \frac{(1-\varphi)}{\xi} + \xi \; \|\nabla \varphi\|^2  \right]  \; d\bfx. \label{reg:energya}
\end{align}
In the above formulation \ref{reg:energya}, $\varphi \in \mbox{H}^{1} \left({\Lambda}; [0, \; 1] \right) $ is introduced as a smooth scalar \textit{phase-field function}. The parameter $\kappa \ll 1$ serves as a numerical regularization parameter specifically for the bulk energy term \cite{heister2015primal}. The constant $c_w$ is the normalization constant. The region where $\varphi \approx 0$ is considered as an approximation of the crack set, while $\varphi=1$ signifies the non-fractured zone. The length parameter $\xi >0$ is responsible for controlling the width of the smooth approximation to the crack set $\Gamma$.

Finally, the quasi-static minimization problem is formulated as follows:

\begin{form1}{form1}
 At an initial time $t=t_0$, with the condition $\varphi( \bfx, \; t_0) =1$ for all $\bfx \in \Lambda$, the solution $\Phi_{\xi}(\bfx, \; t_0)$ is sought as the minimizer:
\begin{equation}
\Phi_{\xi}(\bfx, \; t_0) = \argmin \left\{E_{\xi}(\hat{\Phi}, \; \varphi(\bfx)=1 \; \forall \bfx \in \Lambda) \colon \hat{\Phi} \in H^1 (\Lambda), \; \hat{ \Phi}(\bfx)=g(t_0) \;  \forall \bfx \in \Lambda_D \right\}.
\end{equation}
For subsequent time steps $t=t_k$, where $k=1, \ldots, m$, the pair $\left( \Phi_{\xi}( \bfx, \; t_k), \, \varphi_{\xi}( \bfx, \; t_k) \right)$ is determined by satisfying:
\begin{align}
\left( \Phi_{\xi}(\bfx, \; t_k), \, \varphi(\bfx, \; t_k) \right) = \argmin \big\{ &E_{\xi}(\hat{\Phi},  \, \hat{\varphi} \colon \hat{\Phi}(\bfx) \in H^1 (\Lambda), \; \hat{ \Phi}(\bfx)=g(t_k) \;  \forall \bfx \in \Lambda_D; \notag \\
&\hat{\varphi}(\bfx) \in  H^1 (\Lambda), \;   \partial_t  \hat{\varphi}(\bfx, t) \leq 0 \; \forall \bfx \in \Lambda \big\}. \label{eq_min_con}
\end{align}
\end{form1}
The constraint $\partial_t  \hat{\varphi}(\bfx, t) \leq 0$ in Equation \eqref{eq_min_con} is imposed to enforce crack irreversibility \cite{heister2015primal,manohar2025convergence,manohar2025b}, reflecting the physical reality that cracks cannot heal. Consequently, this formulation results in a time-dependent (or quasi-stationary) variational inequality system. The existence of minimizers for the energy functional $E_{\xi}(\cdot \;, \cdot)$ is contingent upon the functional's convexity. It has been established that the functional exhibits convexity with respect to $\Phi$ when $\varphi$ is held constant, and similarly, with respect to $\varphi$ when $\Phi$ is fixed. In pursuit of understanding this property, the following proposition is put forth:
\begin{proposition}
The energy functional $E_{\xi}(\Phi, \, \varphi)$, defined as:
\begin{equation} \label{eq:prop_energy}
   E_{\xi}(\Phi, \,  \varphi) =  \frac{\mu}{2} \int_{{\Lambda}} \left( (1- \kappa)  \overline{\varphi}^2 + \kappa  \right) \;  \mathcal{W}(\Phi, \, \overline{\varphi}) \; d \bfx + \dfrac{\mathcal{G}_c}{c_w}\int_{{\Lambda}} \left[ \frac{(1-\varphi)}{\xi} + \xi \; \|\nabla \varphi\|^2  \right] \; d\bfx
\end{equation}
is convex when $\overline{\varphi}$ is linearly extrapolated and held fixed.
\end{proposition}
Furthermore, the behavior of the phase-field variable is governed by a crucial maximum principle. This principle is formally stated as follows:
\begin{proposition}[\textbf{The Maximum Principle}]
It is asserted that the phase-field variable $\varphi_{\xi}$ satisfies the maximum principle, such that its values are bounded within a specific range. Specifically, for all spatial points $\bfx \in \Lambda$ and for every discrete time step $k \geq 1$, the following condition holds:
\begin{equation} \label{eq:max_principle}
0 \leq  \varphi_{\xi}( \bfx, \, t_k) \leq \varphi_{\xi}(\bfx, \, t_{k-1}).
\end{equation}
This inequality implies that the phase-field variable, representing the extent of damage or fracture, is always non-negative and monotonically non-increasing over time, reflecting the physical constraint that cracks, once formed, cannot heal or reduce in extent.
\end{proposition}
\begin{remark}
It should be noted that \textbf{Formulation-1} is inherently nonlinear, a characteristic arising from both the nonlinear strain energy density term and the presence of an inequality constraint. The formidable challenge posed by the nonlinear bulk energy term necessitates the design of highly efficient algorithms. Furthermore, the lack of regularity in time, where solutions can exhibit abrupt jumps between successive time steps under varying loads, necessitates the application of specific extrapolation strategies. These strategies are crucial for counteracting the non-convexity of the energy functional, a property extensively discussed in the literature \cite{bourdin2000numerical,bourdin2008variational,bourdin2007numerical,heister2015primal,wick2017modified,fan2021quasi}.
\end{remark}

\subsubsection{Crack irreversibility}
The physical requirement of crack irreversibility in the quasi-static evolution, denoted by $t \mapsto \Gamma(t)$, dictates that the crack set must be non-decreasing with time. This implies that for any $0 \leq t_1 \leq t_2 \leq T$, the condition $\Gamma(t_1) \subseteq \Gamma(t_2)$ must be satisfied.

This aforementioned condition can alternatively be expressed through the phase-field variable as:
\begin{equation}\label{crack-irr-1}
\partial_t \varphi(\bfx, t) \leq 0,
\end{equation}
which explicitly enforces \textit{crack irreversibility}. However, an alternative approach involves imposing an equality constraint. At a given time $t = t_{k-1}$, for $k \in \{1, 2, \ldots, N\}$, a set $\Xi(t_{k-1})$ is defined as:
\begin{equation}
\Xi(t_{k-1}) := \left\{ x \in \overline{\Lambda} \; \colon \; \varphi_{\xi} (\bfx, \, t_{k-1}) < \varkappa \right\},
\end{equation}
where $\varkappa$ is a positive parameter ($0 < \varkappa \ll 1$) whose value is dependent on the regularized length scale $\xi$. It is clear that the set $\Xi(t_{k})$ is non-empty. The crack irreversibility condition, as represented by Equation \eqref{crack-irr-1}, can then be reformulated as:
\begin{equation}\label{crack-irr-2}
\varphi_{\xi}(\bfx, \, t_i) =0 \quad \forall \; \bfx \in \Xi(t_{k-1}) \; \mbox{and} \; \forall \; i \; \mbox{such that} \; k \leq i \leq N.
\end{equation}
This specific reformulation has the advantage of simplifying the minimization problem significantly, as it allows for the setting of $\varphi_{\xi}(\bfx, \, t_i)$ to zero within the identified crack regions, which consequently linearizes the entire process.

\subsection{A three-field formulation for quasi-static crack propagation in strain-limiting solids}
In this quasi-static evolutionary framework, the {crack or damaged region} is conceptualized as a narrow zone where the phase-field, denoted $\varphi_{\xi}(t_k)$, approaches zero. Conversely, the {pristine, undamaged material} is characterized by regions where $\varphi_{\xi}(t_k)$ is close to one. The {regularization parameter $\xi$ dictates the width of this transitional layer}. As $\xi$ tends towards zero, the phase-field approximation is anticipated to become sharper, ideally converging to a discrete crack representation. In this refined state, $v$ would be nearly unity across the undamaged domain and precisely zero on the lower-dimensional crack set.

A significant hurdle in the existing quasi-static phase-field formulations lies in the optimal selection of the regularization length parameter, $\xi$. The current literature often lacks a definitive consensus on its determination. To address this, Phansalkar et al.~\cite{phansalkar2022spatially, phansalkar2023extension} pioneered {adaptive models} where $\xi$ is no longer a fixed constant but treated as an additional, variable field. In their approach, the total energy functional is minimized not only with respect to displacement and the phase-field but also with respect to this spatially varying regularization length. Building upon these adaptive concepts, the present study extends their applicability by incorporating a linear material degradation function, typical of \textsf{AT1}-type models. Specifically, it investigates its behavior under {anti-plane shear loading}. A crucial consideration for such adaptive models is the boundedness of the regularization length field, $\xi$. Consequently, this work introduces a {revised energy functional} that integrates two positive parameters. These parameters are specifically designed to ensure that $\xi$ remains within predetermined upper and lower bounds. The resulting modified energy functional is then formulated as:

\begin{align}\label{reg:energy}
E_{\xi}(\Phi, \,  \varphi, \, \xi) =  & \frac{\mu}{2} \int_{{\Lambda}} \left( (1- \kappa)  {\varphi}^2 + \kappa  \right) \;  \mathcal{W}(\Phi, \, \varphi) \; d \bfx + \dfrac{\mathcal{G}_c}{c_w}\int_{{\Lambda}} \left[ \frac{(1-\varphi + \eta)}{\xi} + \xi \; \|\nabla \varphi\|^2  \right] \; d\bfx \notag \\
& + \int_{\Lambda} \delta \,  \xi \, d\bfx.
\end{align}

Our proposed energy formulation introduces parameters $\eta$ and $\delta$ that exhibit an implicit dependency on $\xi$. This design allows for a strategically smaller $\xi$ near the crack-tip, while employing a slightly larger value in other regions. This adaptive strategy offers a significant edge over conventional methods that rely on a globally constant $\xi$, which often necessitate $\xi > h$ within the damage zone. It's well-established that the stability and convergence of solutions to the Francfort-Marigo model hinge critically on the $\xi$ parameter. Consequently, selecting an optimal $\xi$ is paramount for achieving stable solutions for this model. Our proposed energy is minimized with respect to three coupled field variables: the Airy stress function $\Phi(\bfx)$, the phase-field $\varphi(\bfx)$, and the length scale parameter $\xi$. This minimization can be carried out assuming either a universally constant $\xi$ or a spatially varying field.

\subsection{Determination of optimal global parameter $\xi$}
Our initial analysis begins by considering a {spatially uniform (constant) optimal length scale parameter, $\xi$}, applied across the entire domain $\Lambda$. By minimizing the energy functional with respect to the Airy stress function $\Phi(\bfx)$, the phase-field $\varphi(\bfx)$, and this constant parameter $\xi$, we arrive at the following {system of Euler-Lagrange equations}: These equations represent the governing principles that define the equilibrium state under this simplified assumption.
\begin{subequations}
\begin{align}
    \partial_{\Phi} E_\xi &= \dfrac{\mu}{2}\int_{\Lambda} \left( (1-\kappa) \varphi^2 + \kappa \right)  \partial_{\Phi} \mathcal{W}(\Phi, \, \varphi) \cdot \nabla {\psi} \, d\bfx  &&  \text{in } \Lambda \label{eq:Eu_const_xi} \\[10pt]
   \partial_{\varphi} E_\xi &= \dfrac{\mu}{2} \int_{\Lambda} \partial_{\varphi} \left[     \left( (1-\kappa) \varphi^2 + \kappa \right) \, \mathcal{W}(\Phi, \, \varphi) \right] \psi \; d\bfx - \frac{G_c}{c_v} \int_{\Lambda} \frac{\psi}{\xi}\, d\bfx \\
         &+ \frac{2G_c}{c_v} \int_{\Lambda} \xi \nabla \varphi \cdot \nabla \psi \, d\bfx  &&  \text{in } \Lambda \label{eq:Ev_const_xi} \\[10pt]
 \partial_{\xi}    E_{\xi} &= \frac{G_c}{c_v} \int_{\Lambda}\left( \frac{-(1 - \varphi+\eta)}{\xi^2} + \| \nabla \varphi \|^2 \right) \, d\bfx 
            + \int_{\Lambda} \zeta \, d\bfx &&  \text{in } \Lambda \label{eq:Exi_const_xi} \\[10pt]
    \xi &= \sqrt{ \frac{ \frac{G_c}{c_v} \int_{\Lambda} (1-\varphi+\eta)\, d\bfx}
                 {\frac{G_c}{c_v}\int_{\Lambda} \| \nabla \varphi \|^2 d\bfx + \int_{\Lambda} \delta \,d\bfx }} \label{eq:xi_const_xi}
\end{align}
\end{subequations}
In these equations, $\partial_{\Phi} E_\xi$, $\partial_{\varphi} E_\xi$, and $\partial_{\xi}    E_{\xi}$ represent the variations of the energy functional with respect to the fields $\Phi$, $\varphi$, and $\xi$, respectively. These variations are crucial for deriving the governing equations, often representing equilibrium conditions or optimality criteria within the model. The parameter $\mu$ is an elastic modulus, quantifying the material's resistance to elastic deformation, while $k$ is a residual stiffness parameter, accounting for any remaining stiffness. $G_c$ is the critical energy release rate, a fundamental material property in fracture mechanics. Furthermore, $c_v$ is a model constant, and $\eta$ is a regularization parameter used to enhance numerical stability. The variable $\psi \in H^1(\Lambda)$ denotes a test function within the variational framework. These test functions are crucial for transforming differential equations into their weak forms, facilitating numerical simulations using the finite element method. The optimal constant $\xi$ is determined from the condition $\partial_{\xi}    E_{\xi}=0$, which typically implies a state of equilibrium or an optimal configuration for $\xi$.

\subsection{A spatially dependent (or an inhomogeneous) $\xi$}
To enhance the model's adaptability, we subsequently extend its framework by treating $\xi$ not as a constant, but as a spatially varying field variable, denoted as $\xi(\bfx)$. This means that the value of $\xi$ can now change from one point to another within the material, allowing for the representation of inhomogeneous properties or localized phenomena. This extension necessitates a minimization of the modified energy functional, now expressed with respect to the spatially dependent fields $\Phi(\bfx)$, $\varphi(\bfx)$, and the newly introduced $\xi(\bfx)$. This variational approach yields a set of governing equations: specifically, weak-form Euler-Lagrange equations for $\Phi(\bfx)$ and $\varphi(\bfx)$. These equations are typically partial differential equations that describe the equilibrium or evolution of the primary fields. In contrast, for $\xi(\bfx)$, we obtain a pointwise algebraic optimality condition. This condition, derived from setting the variation of the energy with respect to $\xi(\bfx)$ to zero, means that $\xi(\bfx)$ can be determined locally at each point in space without needing to solve a differential equation for it. This distinction is often beneficial for computational efficiency.

\begin{subequations}
\begin{align}
    \partial_{\Phi} E_\xi &= \dfrac{\mu}{2}\int_{\Lambda} \left( (1-\kappa) \varphi^2 + \kappa \right)  \partial_{\Phi} \mathcal{W}(\Phi, \, \varphi) \cdot \nabla {\psi} \, d\bfx  &&  \text{in } \Lambda \label{eq:Eu_const_xi} \\[10pt]
   \partial_{\varphi} E_\xi &= \dfrac{\mu}{2} \int_{\Lambda} \partial_{\varphi} \left[     \left( (1-\kappa) \varphi^2 + \kappa \right) \, \mathcal{W}(\Phi, \, \varphi) \right] \psi \; d\bfx - \frac{G_c}{c_v} \int_{\Lambda} \frac{\psi}{\xi}\, d\bfx \\
         &+ \frac{2G_c}{c_v} \int_{\Lambda} \xi \nabla \varphi \cdot \nabla \psi \, d\bfx  &&  \text{in } \Lambda \label{eq:Ev_const_xi} \\[10pt]
 \partial_{\xi}    E_{\xi} &= \int_{\Lambda} \left[  \frac{G_c}{c_v}  \left( \frac{-(1 - \varphi+\eta)}{\xi(\bfx)^2} + \| \nabla \varphi \|^2 \right)   +  \delta \right] \, d\bfx &&  \text{in } \Lambda \label{eq:Exi_field_xi} \\[10pt]
    \xi(\bfx) &= \sqrt{ \frac{ \frac{G_c}{c_v}  (1-\varphi+\eta)}
                 {\frac{G_c}{c_v} \| \nabla \varphi \|^2  +  \delta  }} \label{eq:xi_field_xi}
\end{align}
\end{subequations}
In this refined formulation, the governing equations for $\partial_{\Phi} E_\xi$ and $\partial_{\varphi} E_\xi$ maintain their integral forms. These integral equations are typically derived through variational principles, representing the global equilibrium or balance laws for the primary fields $\Phi$ and $\varphi$ across the entire domain. Their integral nature often makes them suitable for numerical solution methods, such as the finite element method, as they can directly incorporate boundary conditions and material inhomogeneities. Crucially, the condition $\partial_{\xi} E_{\xi}=0$, explicitly given by Equation~\eqref{eq:Exi_field_xi}, transforms into a local algebraic equation for $\xi(\bfx)$. Unlike the integral forms for $\partial_{\Phi} E_\xi$ and $\partial_{\varphi} E_\xi$, this condition provides a direct, point-wise relationship for determining the value of $\xi$ at any given spatial location. This means that $\xi(\bfx)$ doesn't require solving a complex differential equation across the entire domain; instead, its value can be computed algebraically at each point, significantly simplifying the computational process. 

The solution for $\xi(\bfx)$ is precisely provided by Equation~\eqref{eq:xi_field_xi}. This derivation assumes that the parameter $\delta$, which is part of the model's constitutive laws, can also be a {spatially varying function}, denoted as $\delta(\bfx)$.  Allowing the parameter $\delta$ to vary spatially, meaning $\delta(\bfx)$, would significantly enhance the model's capacity to represent heterogeneous materials or conditions where material properties change across the domain. This spatial dependence could capture nuanced behaviors, such as varying material stiffness, damage thresholds, or other characteristics that are not uniform throughout the body. For instance, in composites, different plies might exhibit distinct $\delta$ values, or in geological formations, material properties could gradually transition over distance. 
\begin{remark}
While the introduction of a spatially varying $\delta$ presents a compelling avenue for future research, offering a richer and more accurate description of complex material systems, for the current study, we have considered $\delta$ as a global constant. This simplification is often a necessary initial step to establish the model's fundamental behavior before incorporating additional complexities. Nevertheless, exploring $\delta$ as a spatially varying parameter remains an interesting and important direction for future investigations.
\end{remark} 
\section{Numerical discretization via the continuous Galerkin finite element method}\label{fem}

To enable the numerical solution of the continuous governing equations derived in the previous section, the computational domain $\Lambda$ is systematically transformed into a discrete representation. This is achieved by subdividing $\Lambda$ into a finite collection of non-overlapping simplicial elements, collectively forming a tessellation known as the mesh, denoted by $\Th$. This tessellation provides a complete covering of the domain closure, $\overline{\Lambda}$, ensuring that the mesh boundary precisely conforms to the physical domain boundary, $\partial\Lambda$. A critical requirement for this mesh is the conformity condition, which stipulates that the intersection between any two distinct elements, $\tau_i$ and $\tau_j \in \Th$, must be either an empty set, a single shared vertex, or a common edge. This condition is fundamental for ensuring the continuity of basis functions and the accurate assembly of global system matrices.

The geometric integrity of the mesh is paramount for ensuring the accuracy, stability, and convergence properties of the numerical scheme. The size of an individual element $\tau \in \Th$ is quantitatively described by its diameter, $h_\tau = \text{diam}(\tau)$. The overall resolution of the mesh is then characterized by the global mesh size parameter, $h = \max_{\tau \in \Th} h_\tau$. To guarantee numerical well-posedness and to avoid degenerate elements, the entire mesh family is rigorously assumed to be shape-regular. This condition ensures the existence of a uniform positive constant $c_{\varrho}$ such that, for every element within the mesh, the ratio of its diameter to the diameter of its largest inscribed circle remains bounded, i.e., $h_\tau / \varrho_\tau \leq c_{\varrho}$. This geometric constraint is vital for bounding interpolation errors and ensuring the stability of numerical operators.

For standardization and computational efficiency, all local computations are performed on a canonical reference element, $\refface$. Each physical element $\tau \in \Th$ in the computational mesh is uniquely mapped to this normalized reference geometry through an invertible affine transformation, $M_{\tau}: \refface \to \tau$. The chosen reference simplex is the standard unit simplex:
$$
\refface = \{ \refcoord \in \mathbb{R}^d : \hat{z}_k > 0 \text{ for } k=1,\dots,d, \text{ and } \sum_{k=1}^d \hat{z}_k < 1 \}.
$$
This transformation simplifies integration and differentiation operations, allowing them to be performed on a fixed, well-defined geometry.

Upon this established mesh framework, we construct a finite-dimensional function space. Associated with the set of vertices (or nodes) $\{\nodecoord_i\}_{i \in \nodes}$ of the mesh $\Th$, we define a set of nodal basis functions, $\{\basis_i\}_{i \in \nodes}$. These functions possess crucial properties: they are continuous across the entire domain closure $\overline{\Lambda}$, piecewise linear on each element $\tau$, and satisfy the essential Kronecker delta property at the nodes, $\basis_i(\nodecoord_j) = \delta_{ij}$. This property confirms their role as a Lagrange basis, providing a convenient representation for interpolating field variables.

Building upon this foundation, we define the following finite-dimensional approximation spaces at a discrete time $t=t_n$:
\begin{align*}
\FESpace &:= \text{span} \{ \basis_i : i \in \nodes \} = \left\{ \phi_h = \sum_{i \in \nodes} \Phi_i \basis_i : \Phi_i \in \mathbb{R} \right\}, \\
\FESpaceHomD &:= \left\{ \phi_h \in \FESpace : \Phi_i = 0 \text{ for all } i \in \nodesDboundary \right\}, \\
\FESpaceInhomD &:= \left\{ \phi_h \in \FESpace : \Phi_i = g(t_n, \nodecoord_i) \text{ for all } i \in \nodesDboundary \right\},
\end{align*}
where $\nodesDboundary = \{i \in \nodes : \nodecoord_i \in \partial \Lambda \}$ represents the index set of nodes lying on the Dirichlet boundary, and $g(t_n, \nodecoord_i)$ signifies the prescribed non-homogeneous Dirichlet boundary data.

To effectively manage and track evolving crack patterns within the phase-field framework, we introduce a discrete crack set derived from the phase-field variable $\phasefield \in \FESpace$ from the preceding time step $t_{n-1}$ and a user-defined tolerance parameter $\tolCR$. The set of faces considered "cracked" at $t_{n-1}$ is formally defined as:
$$ \FacesCR(t_{n-1}) := \left\{ f \in \Fh : \phasefield(\nodecoord, t_{n-1}) \leq \tolCR \text{ for all } \nodecoord \in \bar{f} \right\}, $$
where $\Fh$ denotes the set of all faces in the mesh. The discrete crack region itself is then established as the union of the closures of these identified cracked faces: $\CrackRegion(t_{n-1}) := \bigcup_{f \in \FacesCR(t_{n-1})} \bar{f}$. This precise definition allows for the construction of a finite element space that intrinsically incorporates the irreversibility of crack growth from the previous time step:
$$ \FESpaceCrack := \left\{ \phi_h \in \FESpace : \phi_h(\nodecoord) = 0 \text{ for all } \nodecoord \in \CrackRegion(t_{n-1}) \right\}. $$
For conciseness and when the specific time $t_n$ is clear from the context, the notations $\FESpaceCrack$ and $\FESpaceInhomD$ may be simplified to $\FESpaceCrackSimple$ and $\FESpaceInhomDSimple$, respectively.
\subsection{Disrete weak formulation} \label{sec:discrete_formulation}
The continuous weak formulation derived in the last section is derived from a three-field Lagrangian approach for quasi-static crack evolution. It serves as the basis for our finite elements-based numerical method. This formulation, an $L^1$-formulation, is built upon the stationarity conditions of an energy functional that couples the elastic energy of the deforming body with the dissipative energy of crack formation. Crack formation is regularized by a phase-field variable, $\varphi$, which smoothly transitions from an intact state ($\varphi=0$) to a fully cracked state ($\varphi=1$). A key aspect of this formulation is the enforced crack irreversibility, which prevents the healing or reversal of accumulated damage represented by the phase-field.

For numerical implementation using the Finite Element Method (FEM), we discretize this continuous formulation. We seek a solution within a finite-dimensional function space  $\FESpace \subset H^1(\Lambda)$, where $h$ represents the mesh size. The solution simultaneously defines the discrete approximations of three field variables: the Airy stress function $\Phi_h$, the phase-field variable $\varphi_h$, and the characteristic length scale parameter $\xi_h$. $\xi_h$ can be a discrete global variable or a spatially varying discrete variable, depending on the model's configuration.

To effectively address the inherent nonlinearities within the governing partial differential equations, a \textit{Picard's iterative scheme} is employed. This method linearizes the problem at each step by leveraging the solution from the previous iteration in the nonlinear terms, thereby transforming the complex nonlinear system into a sequence of more manageable linear problems that can be solved iteratively until a desired convergence criterion is met.

At a fixed time $t=t_n$, and for $m=1,1, \, 2, \, \ldots$ which denotes the Picard's iteration number, the follwing problem seeks to approximate the solution for $\Phi_n \in \FESpace$, and $\varphi_n \in \FESpace$ as given in the following fomrulation:
\begin{dwf}
Given $\mu$, $\kappa$, $G_c$, $c_v$, $\eta$, $\delta$, $\xi_{iv}$, and nd a mesh size $h$ are all known and given the boundary conditions for $\Phi(t_n) = g(t_n)$ on $\partial \Lambda_D$,  then for fixed time $t=t_n$, and for $m=1,1, \, 2, \, \ldots$, find $\Phi_n^m \in \FESpaceInhomD$, $\varphi_n^m \in \FESpace$ such that 
\begin{subequations}
\begin{align}
&\sum_{K \in T_h} \int_K  \left( (1-\kappa) \varphi_{h, \, n-1}^2 + \kappa \right)  \partial_{\Phi} \mathcal{W}(\Phi_{h, n}^{m-1}; \, \Phi_{h, n}^{m} \; ,  \; \varphi_{h, \, n-1}) \cdot \nabla {\psi_h} \, d\bfx 
= 0, 
\quad \forall \psi_h \in  \FESpace\label{eq:discrete_u} \\[1ex]
& \dfrac{\mu}{2} \, \sum_{K \in T_h} \int_K   \partial_{\varphi} \left[     \left((1-\kappa) \left(\varphi_{h, \, n}^m\right)^2 + \kappa \right) \, \mathcal{W}(\Phi_{h, \, n}^m, \, \varphi_{h, \, n}^{m-1}) \right] \psi_h \, d\bfx 
+ \frac{2G_c}{c_v} \sum_{K \in T_h} \int_K \xi(\bfx) \nabla \varphi_{h, \, n}^m \cdot \nabla \psi_h \, d\bfx \notag \\
&= \frac{G_c}{c_v} \sum_{K \in T_h} \int_K \frac{\psi_h}{\xi(\bfx)} \, d\bfx, 
\quad \forall \psi_h \in \FESpace \label{eq:discrete_v}
\end{align}
\end{subequations}
\end{dwf}

In the above formulation, the constant $\xi_{iv}$ is the initially chosen length of the damage zone. This value is typically set as a multiple of the characteristic mesh size, $h$. For instance, one can choose $\xi_{iv}$ to be $2h$ or $5h$, with the specific choice often depending on the coarseness of the mesh and the desired balance between computational cost and accuracy in capturing the damage evolution. 

However, the formulation also allows for a more sophisticated approach where $\xi$ is not a fixed global constant but a spatially adaptive parameter, denoted as $\xi(\bfx)$. This means that the characteristic length of the damage zone can vary across the computational domain. This dynamic adjustment of $\xi(\bfx)$ is particularly powerful for resolving intricate crack patterns or handling situations where the damage zone's size naturally changes based on local material properties or stress concentrations. In the case of spatially adaptive $\xi$, one can dynamically change the value of $\xi(\bfx)$ in the assembly of the equation~\ref{eq:discrete_v}. This dynamic adaptation allows the numerical model to concentrate computational effort and refine the representation of the damage zone precisely where it's needed, leading to more accurate and efficient simulations of crack propagation.

After solving the discrete formulation, we can determine the damage zone length $\xi$. This can be either a globally optimal value or a spatially adaptive value, depending on the specific needs of the simulation. 

\section{Discussion and numerical considerations} \label{rd}

This section bridges our theoretical framework with practical application through a specific numerical example. This example leverages the adaptive algorithms detailed previously to simulate crack evolution. We consider an elastic unit square featuring a pre-existing crack along one edge, subjected to antiplane shear loading on its boundaries. Opting for an antiplane shear model in this quasi-static crack problem provides a crucial and simplified testbed for validating our proposed framework. In this particular setup, the displacement vector is notably simplified, possessing only a single non-zero component that is perpendicular to the analysis plane. This simplification allows us to rigorously test our formulation in a scenario with a scalar-valued displacement unknown, thereby isolating the core mathematical and numerical challenges from the added complexities of full vectorial elasticity.

Our investigation specifically employs a three-field formulation utilizing the \textsf{AT1} phase-field model, which consists of a nonlinear strain energy density function. Unlike traditional approaches that represent cracks as sharp discontinuities, the \textsf{AT1} model uses a continuous scalar field (the phase-field) that smoothly transitions between an undamaged state (represented by a value of 0) and a fully fractured state (represented by a value of 1) over a small, characteristic length scale. A significant advantage of the \textsf{AT1} model, especially when compared to its counterpart, the \textsf{AT2} model, lies in its ability to predict a well-defined elastic limit. This means that damage accumulation only commences after the strain energy surpasses a specific, non-zero threshold. This characteristic often provides a more accurate representation of the behavior observed in many brittle materials.

To effectively solve the coupled system involving the Airy stress function ($\Phi$) and phase-field ($\varphi$) variables, we developed and implemented a multi-step, staggered iterative solver. Our quasi-static analysis progresses by applying a constant $\Phi$ increment at each load step, a strategy that ensures stable convergence of the solution. A critical physical constraint that demands numerical enforcement is the irreversibility of the crack---meaning a crack, once formed, cannot heal. Our code handles this by applying a standard technique: The irreversibility condition is implemented as a local inequality constraint on the phase-field variable, $\varphi$. Specifically, for any given time step $t_n$, the phase-field value $\varphi(t_n)$ is enforced to be less than or equal to its value at the preceding step, $\varphi(t_{n-1})$. This mechanism effectively prevents any increase in the accumulated damage, thereby guaranteeing the irreversible nature of the fracture process. Within the staggered iterative solver, after a potential new phase-field solution is computed for the current step, a straightforward projection is applied. If, at any point, the computed phase-field value indicates a lower damage state than that from the previous converged step, its value is immediately reset to the previous, higher (or equal) value before proceeding to the next iteration. Furthermore, if the calculated damage value falls below a prescribed small threshold, it is set to zero for all subsequent time steps. While this last step could theoretically lead to a crack-widening effect, we did not observe such behavior in our simulations. To further enhance both computational efficiency and solution accuracy, this solver is synergistically combined with a sophisticated adaptive mesh refinement (AMR) strategy. This AMR scheme utilizes a multi-step procedure for both marking elements that require finer resolution and coarsening the mesh in regions where less detail is needed. The entire computational framework was developed in \textsf{C++} and built upon the versatile, open-source finite element library, \textsf{deal.II} \cite{arndt2021deal}.
\subsection{Parameter estimation}
The successful application of our proposed formulation depends on the careful selection of several key parameters, which are crucial for ensuring the convergence, accuracy, and overall efficiency of our method. For the simulations, the physical material properties were defined by a shear modulus of $\mu=80.8$ and a critical energy release rate of $G_c =2.7$. Beyond these fundamental physical constants, several model-specific parameters require precise calibration. The parameter $\kappa$, which serves to regularize the estimation of the bulk energy, was intentionally set to a very small value of $10^{-10}$. This deliberate choice is crucial because it effectively prevents the overestimation of surface energy density in the high-gradient region near the crack tip, thereby ensuring the model's physical accuracy and stability.

The model parameters, $\eta$ and $\delta$, were rigorously calibrated through an analysis of the characteristic profile of the phase-field variable, $\varphi$, under varying conditions. Initially, in regions far removed from the crack, where the material is considered fully intact ($\varphi=1$), and near domain boundaries where the homogeneous Neumann condition ($\mathbf{n} \cdot \nabla \varphi = 0$) naturally causes the phase-field gradient $\| \nabla \varphi \|$ to diminish, a fine resolution of the phase-field is generally unnecessary. In these areas, we assign a larger regularization length scale, specifically $\xi = 5h$, directly tying it to the local mesh size $h$. This choice reflects a strategic coarsening of the phase-field representation. When this value is substituted into the model's governing relations, it yields a direct expression for the parameter $\eta$ as a function of the mesh size:
\begin{equation}
    \eta=\frac{100h^2c_v\delta}{G_c} \label{eq:eta_calib}
\end{equation}
In stark contrast, a distinct approach is necessitated at the crack interface itself. At the very tip of the crack, the material is considered fully damaged (i.e., $\varphi=0$), and there exists a sharp transition zone where the phase-field undergoes rapid variation. To accurately resolve this critical high-gradient region, a smaller regularization length that precisely scales with the mesh resolution is indispensable. Based on the observed inverse proportionality between the gradient norm $\| \nabla \varphi \|$ and the mesh size $h$, we define a smaller regularization length: $\xi = 2h$. Substituting this value into the governing equations provides a specific formula for the parameter $\alpha$:
\begin{equation}
    \delta=\frac{3G_c}{96c_v h^2} \label{eq:alpha_calib}
\end{equation}
Collectively, these meticulously derived relationships enable the systematic calculation of all required parameters for any given mesh discretization. The resulting values for various mesh sizes, demonstrating their calibration, are concisely summarized in Table \ref{tab:param_calib}.

\begin{table}[h!]
    \centering
    \caption{Calibrated model parameters for different mesh sizes.}
    \label{tab:param_calib}
    \begin{tabular}{cccc}
        \toprule
        {$\#$ of cells} & {h} & {${\delta}$} & {${\eta}$} \\
        \midrule
        128 & 0.008 & 493.75 & 9.36 \\
        256 & 0.004 & 1975 & 9.36 \\
        512 & 0.002 & 7900 & 9.36 \\
        \bottomrule
    \end{tabular}
\end{table} 

\subsection{Staggered algorithm }
To efficiently manage the complex, coupled phase-field evolution under antiplane shear in a specially characterized algebraically nonlinear material, we implemented a robust staggered solution algorithm. This computational strategy skillfully breaks down the comprehensive problem into two distinct sub-problems. These sub-problems are then solved sequentially within each discrete time step or boundary load increment. Initially, the momentum balance equation is solved to determine the out-of-plane displacement field. During this stage, the phase-field variable is held constant, taking its value from the preceding iteration. Here, the phase-field functions as a static degradation function, locally reducing the material's shear stiffness and thereby representing the current state of material damage. Following this, with the newly computed displacement field at hand, the strain energy density is calculated. This energy density then acts as the primary driving force for the second sub-problem: solving the phase-field evolution equation. This subsequent update precisely defines the new topology of the diffuse crack network. This alternating solution procedure continues iteratively until a predefined convergence criterion is met. This approach offers a significant advantage over a monolithic scheme, proving to be both computationally efficient and memory-friendly by avoiding the assembly and inversion of a large, fully coupled Jacobian matrix. To further enhance the efficiency of phase-field fracture simulations, this work incorporates an adaptive iterative scheme. The algorithm initiates with a coarser mesh and proceeds to sequentially solve for the displacement field, the phase-field, and a variable length scale. To ensure high accuracy while simultaneously minimizing computational costs, the mesh is dynamically refined in regions where the solution exhibits high gradients. This adaptive strategy enables a robust and efficient simulation of crack propagation. The entire staggered algorithm, including its AMR capabilities, is meticulously outlined in Algorithm 1.

\begin{algorithm}[H]
\caption{Staggered solution algorithm with adaptive mesh refinement based on the damage length scale variable $\xi$}
\begin{algorithmic}[1]
\State Input: Initial mesh, model parameters, boundary condition increment \(\Delta \Phi\), tolerance \(tol\)
\State Input: AMR parameters \(\xi_{\text{refine}}\), \(h_{\min}\); function \(g(\xi)\); flag \(\texttt{refinement\_enabled}\)
\State Output: Converged fields \(\Phi\), \(\varphi\)
\State Initialize damage field \(\varphi_0=1.0\) and initial guess for Airy stress function is obtained by solving linear problem (i.e. $\beta=0$).
\State Set load step counter \(n = 0\).

\While{load step \(n < N_{\text{max}}\)} \Comment{Main loop over load increments}
    \State \textit{// Initialize for staggered iterations at current load step}
    \State Set staggered iteration counter \(k = 0\).
    \State Set initial guesses: \(\Phi_{n,k} \gets \Phi_{n-1}\), \(\varphi_{n,k} \gets \varphi_{n-1}\).
    \Repeat
        \State \textit{// Store previous iteration's solution for error calculation}
        \State \(\Phi_{\text{prev}} \gets \Phi_{n,k}\), \(\varphi_{\text{prev}} \gets \varphi_{n,k}\)
        \State \(k \gets k + 1\)

        \State Solve for Airy stress function \(\Phi_{n,k}\) using \(\varphi_{\text{prev}}\).
        
        \State Solve for phase-field \(\varphi_{n,k}\) using the newly computed \(\Phi_{n,k}\).

        \State \textit{// Check for convergence of the staggered scheme}
        \State Compute $L_2$ norm of displacement error: \( \text{err}_\Phi = \frac{\| \Phi_{n,k} - \Phi_{\text{prev}} \|_2}{\| \Phi_{n,k} \|_2} \)
        \State Compute $L_2$ norm of phase-field error: \( \text{err}_\varphi = \frac{\| \varphi_{n,k} - \varphi_{\text{prev}} \|_2}{\| \varphi_{n,k} \|_2} \)
        
    \Until{\( (\text{err}_\Phi < tol) \) and \( (\text{err}_\varphi < tol) \)} \Comment{Staggered loop converges}

    \State \textit{// Solution at step n has converged: \(\Phi_n \gets \Phi_{n,k}\), \(\varphi_n \gets \varphi_{n,k}\)}
    
    \State \textit{// AMR based on converged phase-field solution}
    \If{\texttt{refinement\_enabled}}
        \State Set \texttt{refinement\_occurred} = False
        \For{each element \(e\) in the mesh}
            \State Let \(\xi_e\) be the local value computed in element \(e\).
            \State Let \(h\) be the current size of element \(e\).
            \If{$\xi_e < \xi_{\text{refine}}$ and $h > h_{\min}$}
                \State Flag element \(e\) for refinement.
                \State Set \texttt{refinement\_occurred} = True
            \EndIf
        \EndFor
        \If{\texttt{refinement\_occurred}}
            \State Refine all flagged elements and project fields \(\Phi_n, \varphi_n\) onto the new mesh.
        \EndIf
    \EndIf

    \State \textit{// Prepare for the next load step}
    \State Update boundary conditions with the next load increment, e.g., \(\Phi^{BC}_{n+1} \gets \Phi^{BC}_{n} + \Delta \Phi\).
    \State \(n \gets n + 1\).
\EndWhile
\end{algorithmic}
\end{algorithm}

\subsection{Problem setup and loading conditions}
Our numerical investigations employ a square computational domain, centrally featuring a pre-existing edge crack situated on its top boundary, as depicted in Figure \ref{fig:crackdomain}. To simulate a {Mode III (tearing) fracture} scenario, we subject the material to {anti-plane shear loading}. This is precisely achieved by imposing time-dependent, opposing displacements along the top boundary, specifically on either side of the crack mouth ($\Gamma_3$):
\begin{align*}
\Phi =
\begin{cases}
-c\,t & \text{for } x \in (0, 0.5) \quad \text{and} \quad y=1, \\
~~c\,t & \text{for } x \in (0.5, 1) \quad \text{and} \quad y=1. \\
\end{cases}
\end{align*}
The selection of boundary conditions is critical for accurately reflecting the underlying physical phenomena:
\begin{itemize}
    \item \textbf{Mechanical boundary conditions:} Beyond the actively loaded segment of $\Gamma_3$, all external surfaces, including the existing crack faces ($\Gamma_c$), are designated as {traction-free} ($\mathbf{n}\cdot \nabla \Phi=0$). This implies that no external forces act on these surfaces, allowing them to deform unhindered in response to the internal stress distribution.
    \item \textbf{phase-field boundary conditions:} We enforce a {homogeneous Neumann condition} ($\mathbf{n}\cdot \nabla \varphi=0$) on all domain boundaries for the phase-field variable. This condition is crucial; it ensures that the evolution and propagation of the crack are driven purely by the material's internal energy state, rather than being artificially influenced or constrained by the domain's edges. In a physical sense, these boundaries are considered 'neutral,' meaning they neither promote nor inhibit crack initiation or growth, thereby enabling the crack's path and velocity to be dictated solely by the evolving stress field within the domain.
\end{itemize}
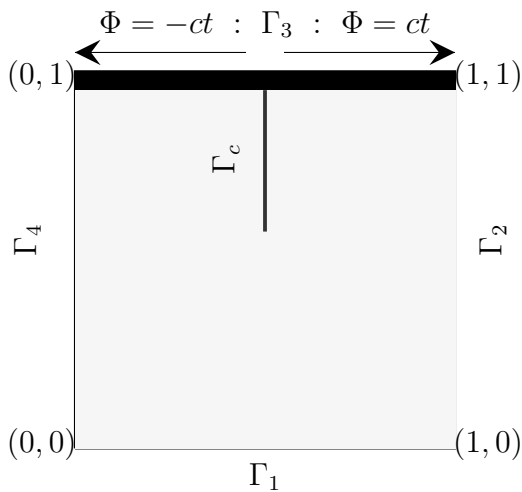
\begin{figure}[H]
\centering
\begin{tikzpicture}[scale=1.25]
    \tikzset{myptr/.style={decoration={markings,mark=at position 1 with %
    {\arrow[scale=3,>=stealth]{>}}},postaction={decorate}}}
    \filldraw[draw=black, thick] (0,0) -- (4,0) -- (4,4) -- (0,4) -- (0,0);
    \shade[inner color=gray!7, outer color=gray!7] (0,0) -- (4,0) -- (4,3.8) -- (0,3.8) -- (0,0);
    \node at (-0.35,0.05)   {$(0,0)$};
    \node at (4.35,0.05)   {$(1,0)$};
    \node at (-0.35,3.95)   {$(0,1)$};
    \node at (4.35, 3.95)   {$(1,1)$}; 
    \node at (-0.5, 2.5)[anchor=east, rotate=90]{$\Gamma_4$};
    \node at (4.4, 2.5)[anchor=east, rotate=90]{$\Gamma_2$};
    \node at (2, -0.3) {$\Gamma_{1}$};
    \node at (2, 4.5) {$ \Phi=-ct ~: ~\Gamma_{3} ~:~ \Phi=ct$};
    \draw [myptr](1.8, 4.2)--(0.0, 4.2);
    \draw [myptr](2.2, 4.2)--(4.0, 4.2);
    \draw [line width=0.5mm, black!80]  (2,3.8) -- (2,2.3);
    \node at (1.59,3.35)[anchor=east, rotate=90]{$\Gamma_c$};
   
\end{tikzpicture}
\caption{A computational domain showing an edge crack with anti-plane shear loading.}\label{fig:crackdomain}
\end{figure}

The remainder of this paper is dedicated to an in-depth investigation of how the phase-field length scale parameter, $\xi$, impacts fracture simulations. Our study is structured around a comparative evaluation of two fundamentally different approaches to managing this parameter. Firstly, we consider a {traditional modeling paradigm} where $\xi$ is maintained as a {constant value throughout the domain}. This global constant is carefully chosen and optimized to ensure an adequate and consistent regularization of the crack interface, a standard practice in many phase-field applications. This framework serves as a baseline against which to evaluate more advanced techniques. Secondly, we introduce and analyze an innovative framework in which $\xi$ is treated as a {spatially varying, inhomogeneous field}. This advanced treatment allows the model to inherently control the characteristic width of the diffuse crack based on local material conditions. By dynamically adjusting its value in response to the computed stress and damage fields, this adaptive approach facilitates a more accurate representation of the physical crack features. The forthcoming analysis of the numerical results will systematically compare these two distinct methodologies. Our objective is to elucidate the significant benefits of the adaptive length scale, particularly in achieving enhanced predictive accuracy for crack propagation and realizing notable gains in computational performance. This rigorous comparison will provide compelling evidence for the efficacy of variable-length scale formulations.

The parameter $\beta$ is of great importance for guaranteeing the numerical algorithm's stability and reliable convergence. To thoroughly investigate its influence, an extensive sensitivity analysis was performed. Our findings indicate that when the value of $\beta$ is greater than $0.001$, the computational code consistently fails to yield a physically meaningful crack propagation path. It is our current understanding that this computational limitation does not originate from the fundamental physical model. Instead, it is primarily attributable to the characteristics of the linear solver presently in use. Currently, a direct solver is employed; however, for future developments, a more formidable and efficient iterative solver, enhanced by an appropriate preconditioner, is planned for implementation to overcome these challenges.
\subsection{Optimal and Globally Constant $\xi$}
This section details the numerical procedure for identifying an optimal and globally constant regularization parameter, $\xi$. The approach involves solving the regularized phase-field model, utilizing a revised version of Algorithm-$1$. A significant modification to the algorithm is the exclusion of local mesh refinement; instead, a uniformly global mesh is maintained throughout the entire computation. For each time step, an optimal value for $\xi$ is computed for the entire domain, determined from the numerical solution as per Equation~\ref{eq:xi_const_xi}.

Regarding the numerical simulations, the material properties were established with a fracture toughness of $G_c = 2.7$ and a shear modulus of $\mu=80.8$~GPa. All other parameters within the model were adopted directly from the calibrated values presented in Table~\ref{tab:param_calib}.

\begin{figure}[H]
    \centering 
    \begin{subfigure}{0.3\textwidth}
        \centering
        \includegraphics[width=\linewidth]{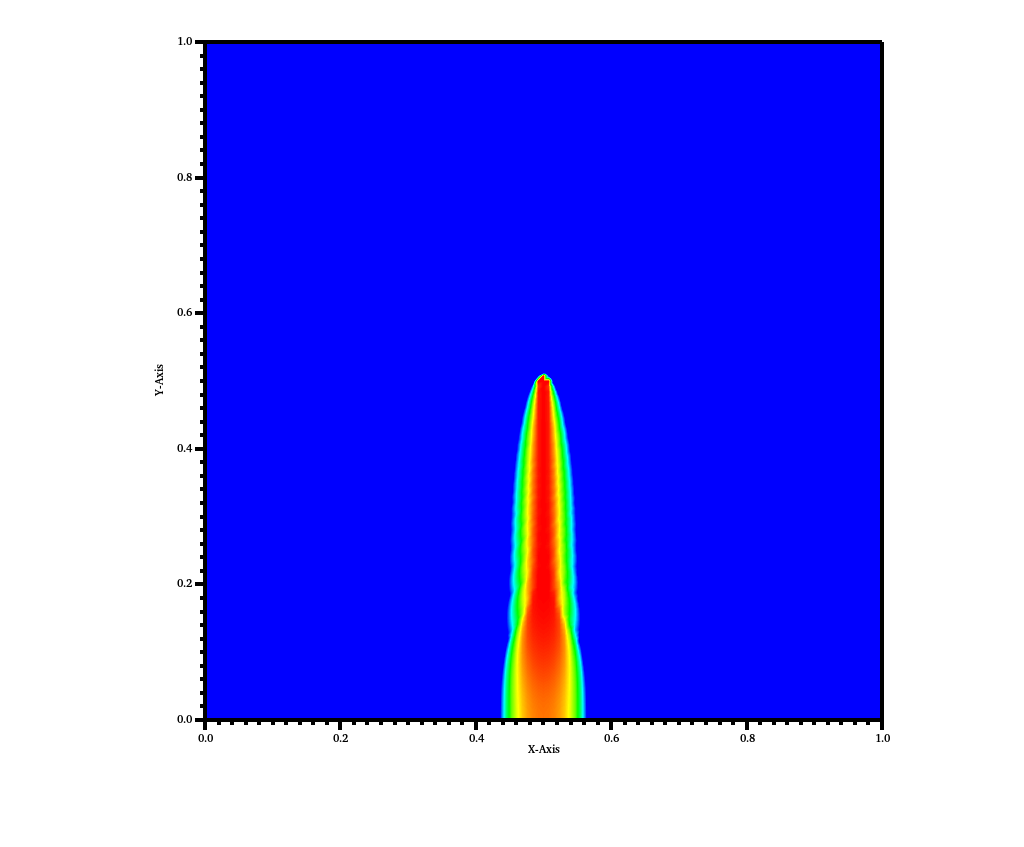}
        \caption{$128$ cells}
        \label{fig:sub1}
    \end{subfigure}
    \hfill 
    \begin{subfigure}{0.3\textwidth}
        \centering
        \includegraphics[width=\linewidth]{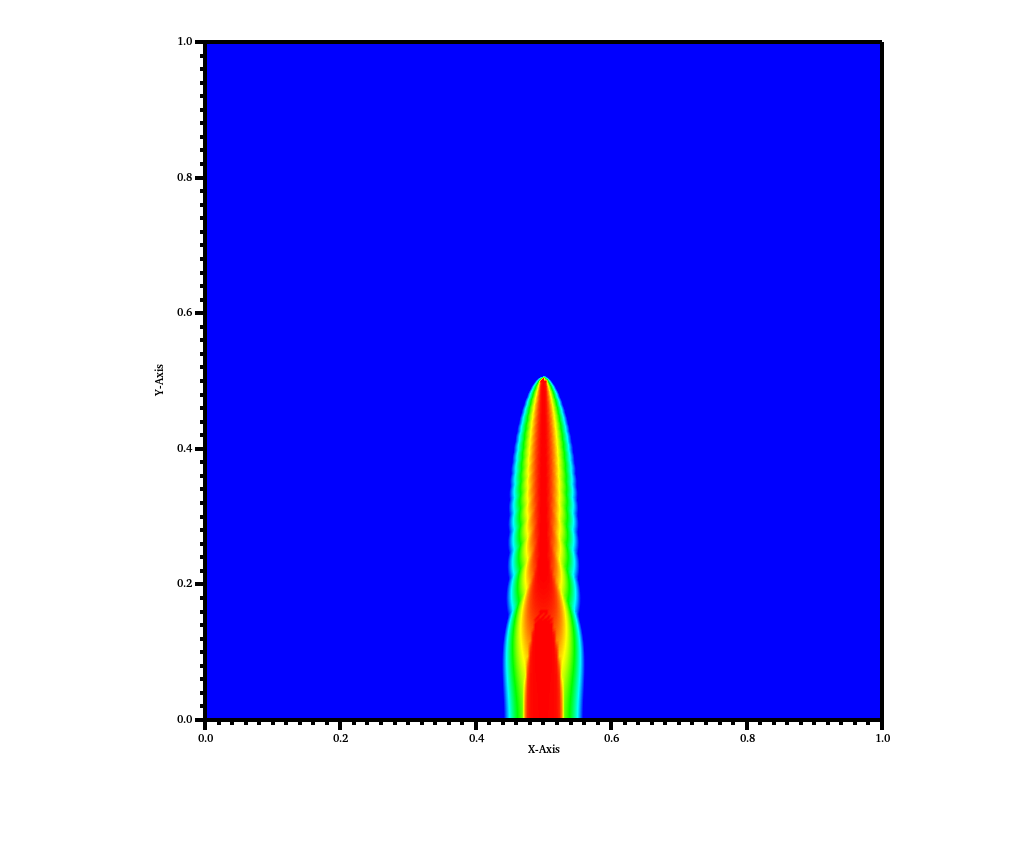}
        \caption{$256$ cells}
        \label{fig:sub2}
    \end{subfigure}
    \hfill 
    \begin{subfigure}{0.3\textwidth}
        \centering
        \includegraphics[width=\linewidth]{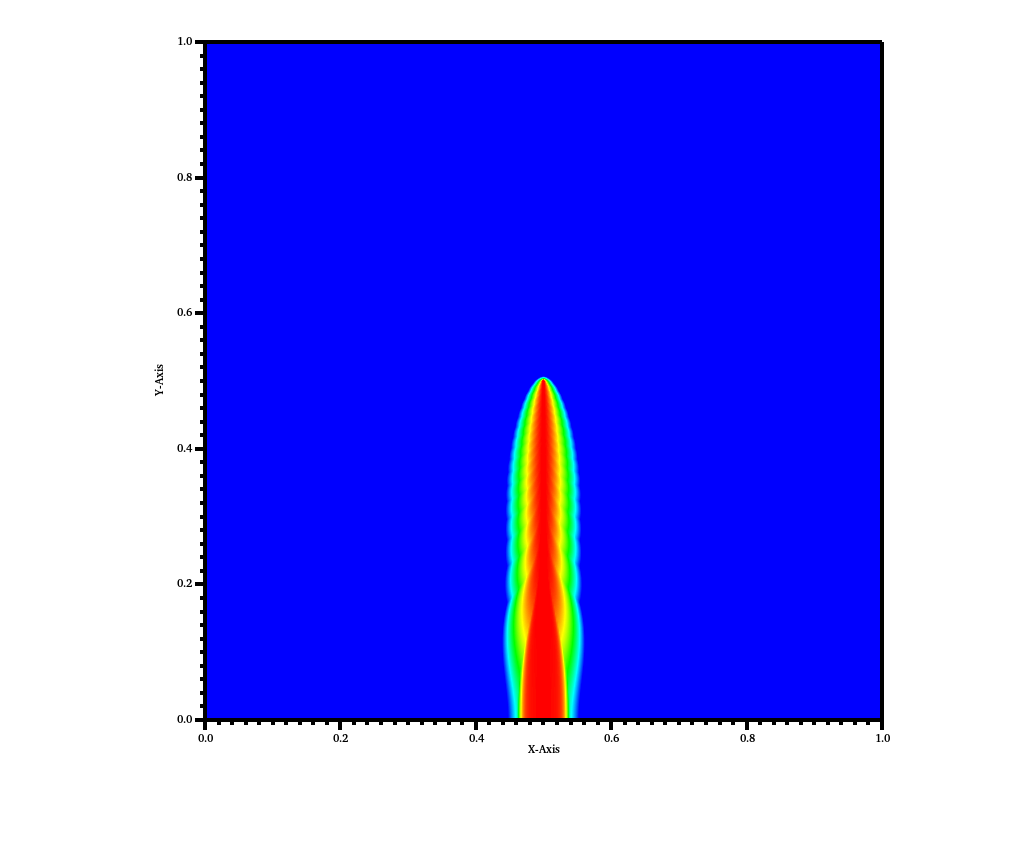}
        \caption{$512$ cells}
        \label{fig:sub3}
    \end{subfigure}
    \caption{Plof of $\varphi$ on the different meshes with globally constant $\xi$.}
    \label{v_globalConsXI}
\end{figure}

In Figure~\ref{v_globalConsXI}, we present a comprehensive visualization of the numerical outcomes, specifically detailing the phase-field variable, $\varphi$. These results were meticulously obtained from simulations conducted across a range of distinct mesh resolutions, allowing for an assessment of solution convergence and detail capture. The contour plots, in particular those for the phase-field variable $v$, serve to graphically illustrate the evolving and predicted crack topology within the simulated material. A clear color scheme is employed within these phase-field plots to delineate material states: regions rendered in red, where the phase-field variable $\varphi$ approximates zero ($\varphi \approx 0$), unequivocally represent areas of completely fractured or damaged material. Conversely, the blue regions, where $v$ approaches one ($\varphi \approx 1$), correspond to the pristine, undamaged state of the material. Concurrently, the accompanying plots for $u$ vividly depict the out-of-plane displacement field, a direct and characteristic response to the anti-plane shear loading conditions that were applied in these simulations.

\begin{table}[h!]
    \centering
    \caption{The optimal value of $\xi$ computed from Equation~\ref{eq:xi_const_xi}.}
    \label{tab:xi}
    \begin{tabular}{cccc}
        \toprule
        {$\#$ of cells } & {h} & {${\xi = m \cdot h}$} & {optimal ${\xi}$} \\
        \midrule
        128 & 0.008 & 0.0390625 &0.13605 \\
        256 & 0.004 & 0.03125   & 0.06927 \\
        512 & 0.002 & 0.0234375   & 0.039844 \\
        \bottomrule
    \end{tabular}
\end{table}

The findings of our comprehensive analysis are concisely presented in Table~\ref{tab:xi}. This table provides a direct comparison between the optimal regularization parameter, $\xi$, meticulously calculated using Equation~\ref{eq:xi_const_xi}, and the comparatively smaller, predetermined value initially employed at the commencement of the minimization procedure. A consistent observation across all evaluated mesh resolutions is that the optimally determined $\xi$ consistently proves to be significantly larger than its pre-selected counterpart. This particular outcome represents a considerable computational advantage, as a larger regularization parameter effectively leads to a relaxation in the overall stiffness of the numerical system, which, in turn, directly translates to a notable acceleration in computational performance and reduced processing times.

To rigorously validate the efficacy and reliability of this proposed approach, supplementary simulations were conducted utilizing the established, albeit smaller, conventional values for the $\xi$ parameter. Crucially, the results from these validation tests demonstrated strikingly similar crack profiles and displacement fields when compared to those obtained with the optimized parameter. This compellingly confirms that the gains in computational efficiency achieved through a larger $\xi$ do not, in any way, detrimentally affect the predictive accuracy or integrity of the simulation results. Consequently, the methodology presented here offers a highly attractive and robust framework for systematically selecting the critical $\xi$ parameter within the \textsf{AT1} model, resulting in substantial reductions in computational time and resource allocation without any discernible compromise in the precision or reliability of the predictions.

\begin{figure}[htb!]
    \centering
    \begin{subfigure}{0.32\textwidth}
        \centering
        \includegraphics[width=\linewidth]{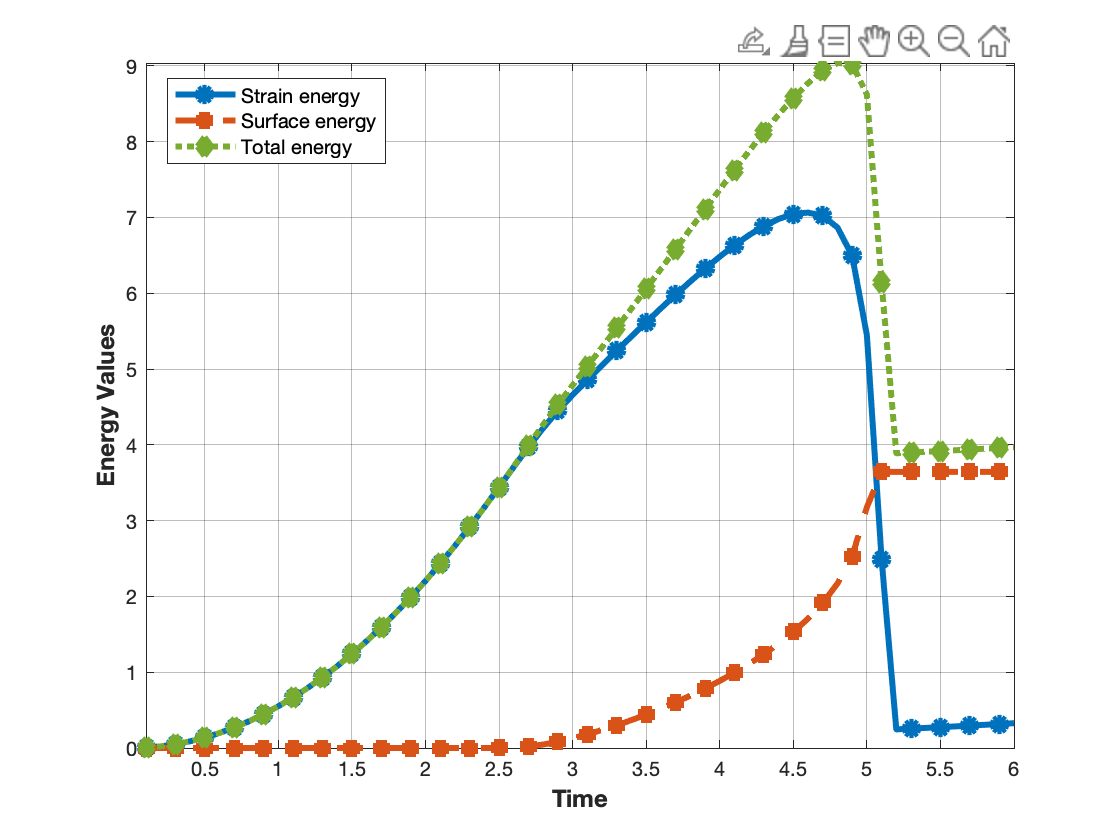}
        \caption{$128$ cells}
        \label{fig:sub1}
    \end{subfigure}\hfill
    \begin{subfigure}{0.32\textwidth}
        \centering
        \includegraphics[width=\linewidth]{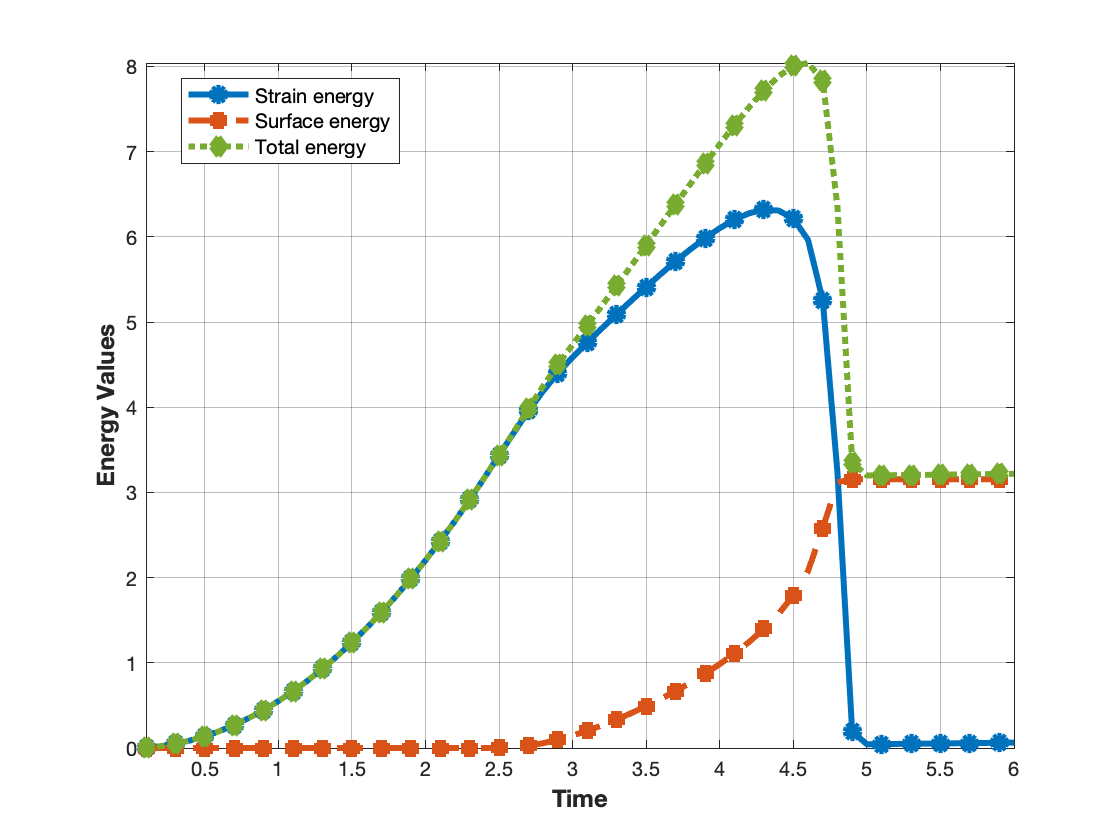}
        \caption{$256$ cells}
        \label{fig:sub2}
    \end{subfigure}\hfill
    \begin{subfigure}{0.32\textwidth}
        \centering
        \includegraphics[width=\linewidth]{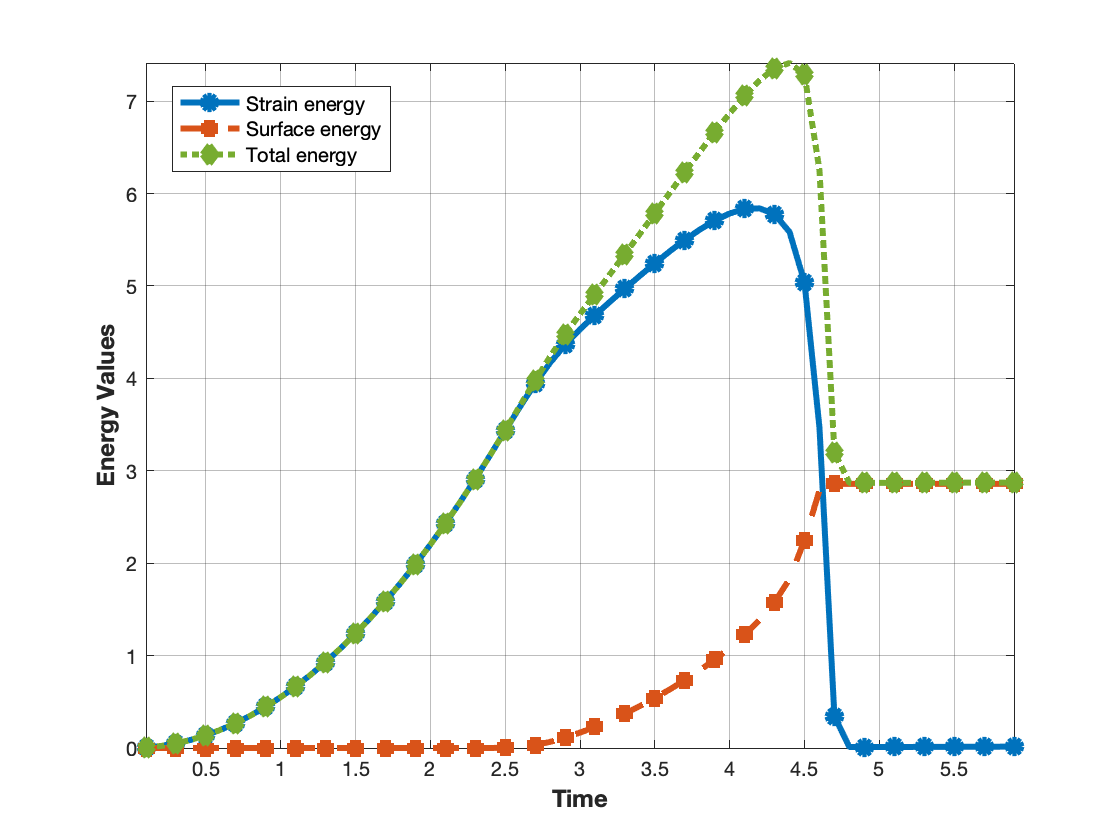}
        \caption{$512$ cells}
        \label{fig:sub3}
    \end{subfigure}
\caption{Evolution of key energy components (elastic strain, newly generated surface, and their sum, the total energy) derived from an \textsf{AT1} phase-field fracture simulation, presented across various computational mesh resolutions.} 
    \label{glob_mesh_energies}
\end{figure}

The plots \ref{glob_mesh_energies} vividly illustrate the fundamental physical process of energy dissipation during fracture propagation: they demonstrate the progressive conversion of stored elastic strain energy within the material into surface energy as the crack initiates and extends. A crucial observation from these results is the clear convergence of the simulated energy profiles, particularly the total energy, as the computational mesh undergoes successive refinement, thereby affirming the numerical robustness and accuracy achieved with finer discretizations.

\subsection{Spatially varying (or inhomogeneous) $\xi$}
This section focuses on solving the minimization problem for the coupled fields of Airy stress function ($\Phi$), the phase-field ($\varphi$), and, crucially, a spatially varying (or inhomogeneous) damage length scale parameter, $\xi$. Our computational approach leverages Algorithm~1, initially implemented on a relatively coarse base mesh comprising $64$ cells. Subsequently, an AMR strategy is integrated, which dynamically adjusts $\xi$ based on local solution characteristics. The AMR process is constrained by an upper bound, allowing for a maximum of four levels of refinement for any individual cell flagged for increased resolution.

This controlled refinement strategy establishes distinct bounds for the effective mesh size, $h$, ranging from a coarsest resolution of $1/2^6 \leq h$ to a finest resolution of $h \leq 1/2^{10}$. These mesh bounds, in turn, induce corresponding constraints on the dynamically varying parameter $\xi$, which is observed to range approximately from $0.15$ down to $0.011$. A critical feature of our AMR methodology is its inherent ability to prevent unbounded mesh growth; it intelligently incorporates a coarsening mechanism for cells that might otherwise become excessively fine. Specifically, cells flagged for refinement beyond the smallest prescribed mesh size are coarsened, thereby ensuring that the AMR process concludes within a finite number of refinement steps, preventing infinite loops and maintaining computational efficiency.

In terms of material properties for the numerical simulation, a fracture toughness of $G_c = 2.7$ and a shear modulus of $\mu=80.8$  were previously defined. The other requisite parameters for the model were configured as follows: $\eta = 31640.6$, $\delta = 3.125$, and $\kappa = 1.0 \times 10^{-10}$. The simulation progressed with a constant time step of $0.01$. It was executed until the phase-field variable, $\varphi$, which serves as a representation of the propagating crack, fully traversed the computational domain. A particularly significant element of this investigation is the role of the length scale parameter, $\xi$. Traditionally, the selection of $\xi$ in phase-field models is guided by its direct proportionality to the local mesh size, $h$. For example, beginning with an initial coarse mesh size of $h=0.015625$, conventional practice would suggest a corresponding $\xi$ value of $10h = 0.15625$. Conversely, for the finest mesh resolution employed in our study ($h=0.00097656$), this standard rule would necessitate a considerably smaller $\xi$ value of $0.0097656$.

In stark contrast to this conventional approach, our simulation dynamically determined and optimized the value of $\xi$ throughout the computation. At the very first time step of the simulation, the dynamically calculated optimal $\xi$ was determined to be approximately $0.03464$. Subsequently, as the crack initiated and began to propagate (evidenced by the decrease of the phase-field variable $\varphi$ from its initial undamaged value of 1), the optimal $\xi$ converged and stabilized within a narrow range, typically varying between $0.02419$ and $0.03644$. It is of considerable importance to emphasize that these computationally derived optimal values are markedly larger than the $0.0097656$ value that would typically be applied when using such a fine mesh in other related studies. This observation represents a significant and key departure from conventional phase-field modeling practice, as further discussed in related literature \cite{yoon2021quasi}.

\begin{figure}[H]
    \centering
    \begin{subfigure}[b]{0.48\textwidth}
        \centering
        \includegraphics[width=\textwidth]{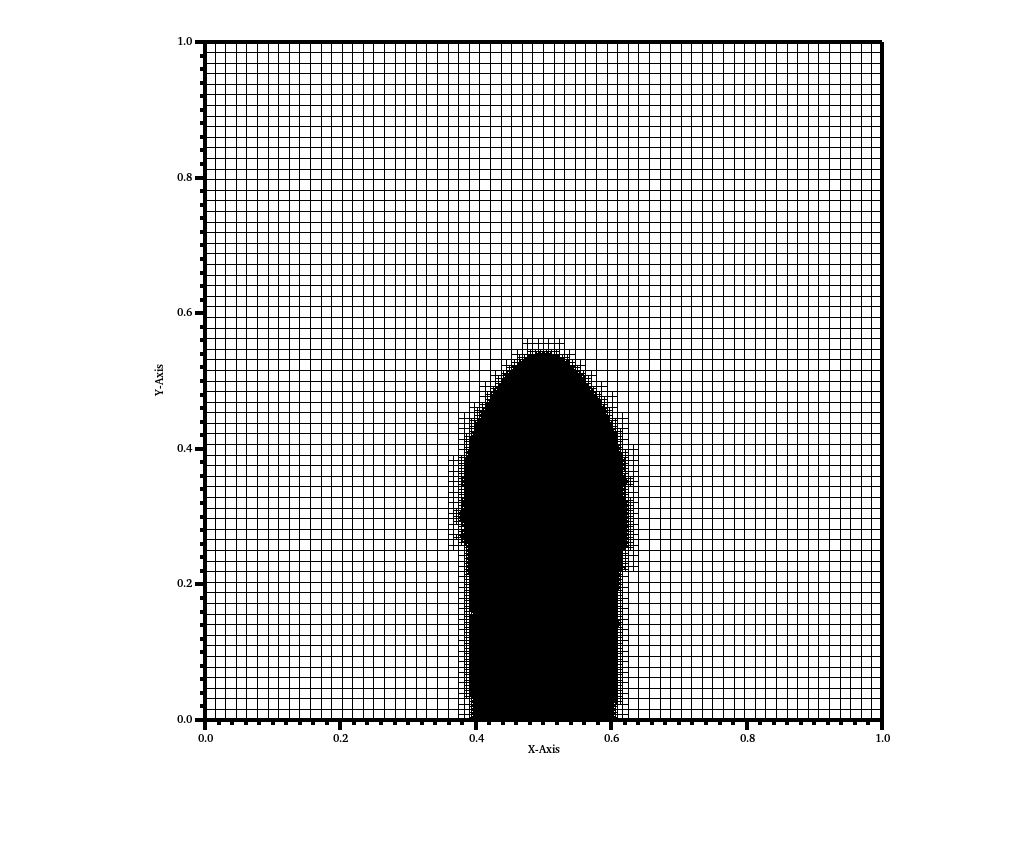}
        \caption{Mesh}
        \label{fig:mesh}
    \end{subfigure}
    \hfill 
    \begin{subfigure}[b]{0.48\textwidth}
        \centering
        \includegraphics[width=\textwidth]{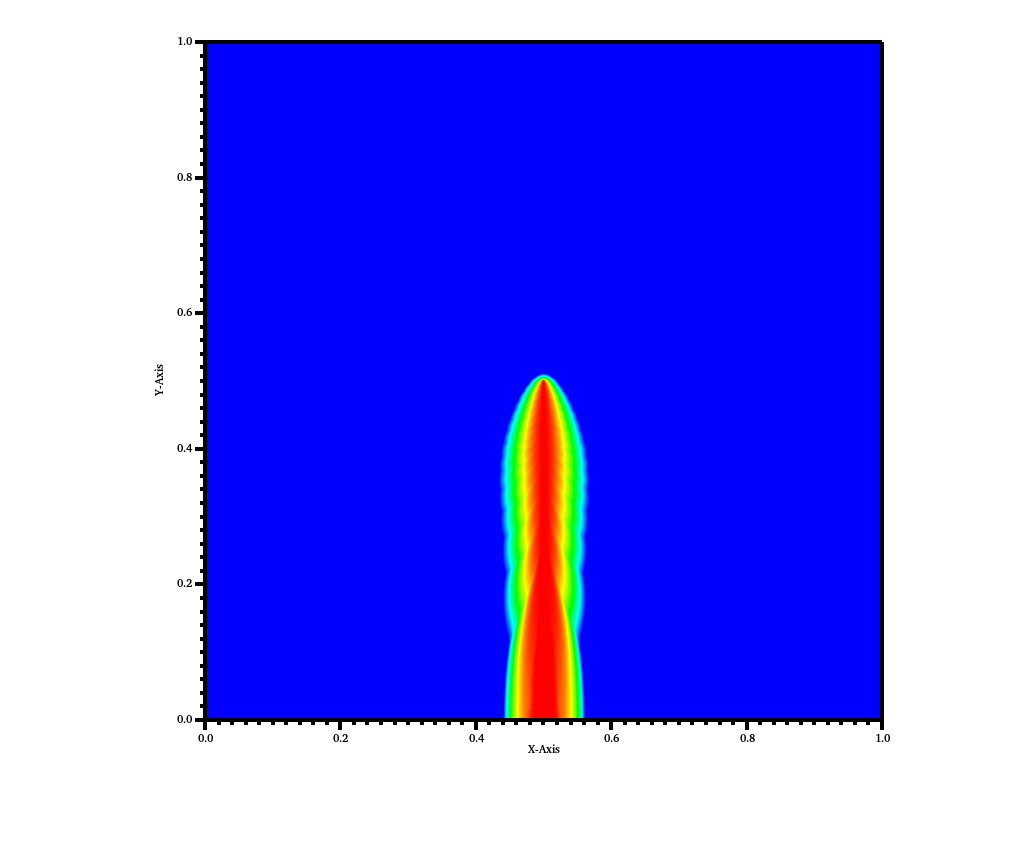}
        \caption{phase-field $v$}
        \label{fig:phase_field}
    \end{subfigure}
\caption{Computational mesh and the phase-field, $\varphi$. These simulation results correspond to a locally varying $\xi$.}
    \label{fig:mesh_v}
\end{figure}

Figure~\ref{fig:mesh_v} illustrates the simulation results obtained using a locally varying length-scale parameter, $\xi$. The computational mesh, shown on the left, was adaptively refined based on the $\xi$ values calculated using Equation~\eqref{eq:xi_field_xi}. This strategy ensures high resolution where the phase-field $\varphi < 1$, leading to a significantly more accurate solution than what a uniformly coarse mesh would provide. The right side of the figure displays the final phase-field, $\varphi$, which visualizes the fully developed crack path: dark red ($\varphi=0$) represents completely fractured material, and blue ($\varphi=1$) indicates undamaged material.

\begin{figure}[H]
    \centering
    \begin{subfigure}[b]{0.48\textwidth}
        \centering
        \includegraphics[width=\textwidth]{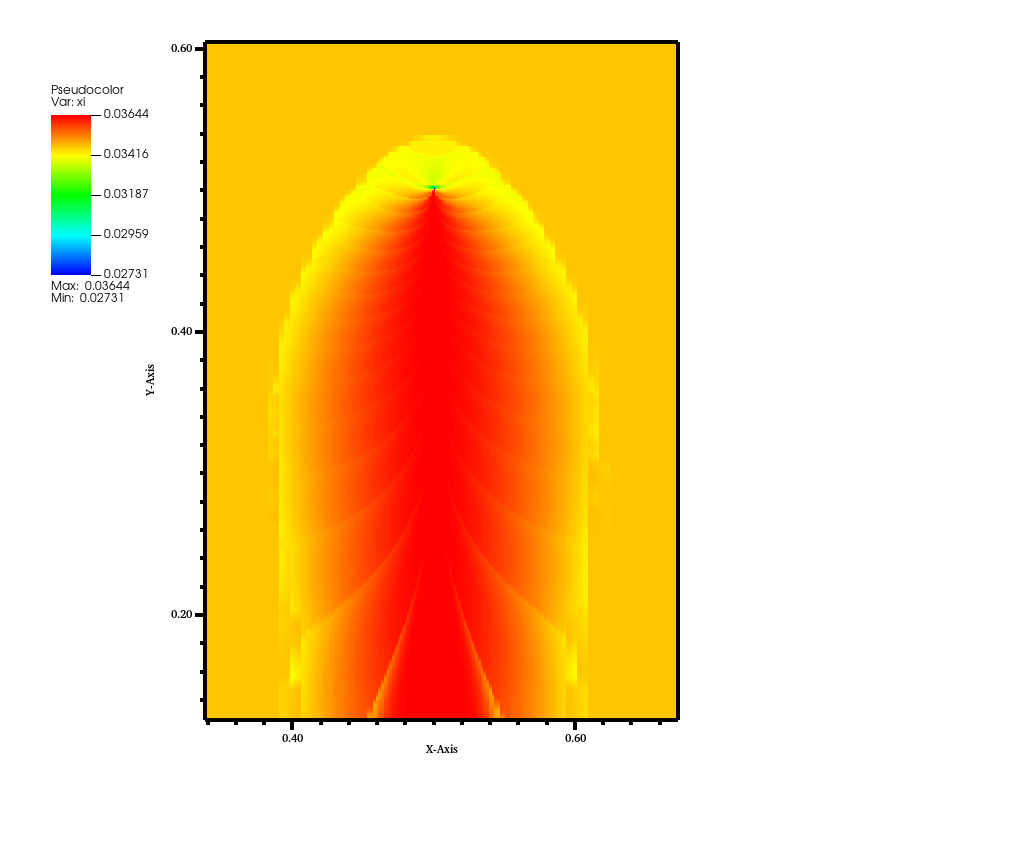}
        \caption{$\xi$ }
        \label{fig:mesh}
    \end{subfigure}
    \hfill 
    \begin{subfigure}[b]{0.48\textwidth}
        \centering
        \includegraphics[width=\textwidth]{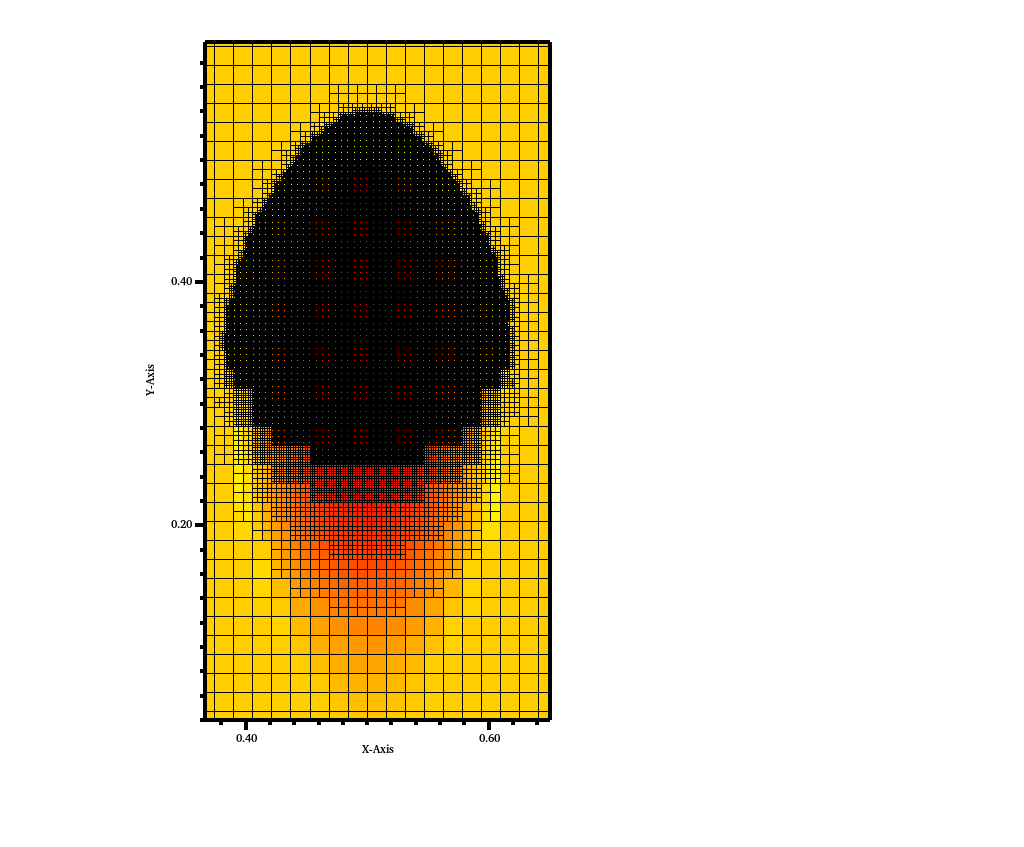}
        \caption{$\xi$ close to initial crack-tip}
        \label{fig:phase_field}
    \end{subfigure}
\caption{ The simulation results correspond to a locally varying $\xi$.}
    \label{fig:loc_xi}
\end{figure}

Figure~\ref{fig:loc_xi} illustrates the dynamic evolution of the locally varying parameter $\xi$. Initially, $\xi$ remains constant at $0.03464$ at the start of the simulation, and this value subtly increases as $\varphi$ decreases. This initial value significantly exceeds the $10h_f$ threshold, where $h_f$ denotes the finest mesh size for a 1024-cell grid. As the variable $v$ begins to decrease from its initial value of $1$, $\xi$ slightly increases away from the initial crack tip, while simultaneously decreasing in the region surrounding the crack. Nearing the final simulation step, $\xi$ reaches its minimum value of $0.02419$. Notably, this minimum is still an order of magnitude larger than typical globally constant $\xi$ values. This demonstrates a key advantage of employing a locally varying $\xi$ within the \textsf{AT1} model, as it facilitates more adaptive and physically accurate parameterization.

\begin{figure}[H]
    \centering
    \includegraphics[width=0.6\textwidth]{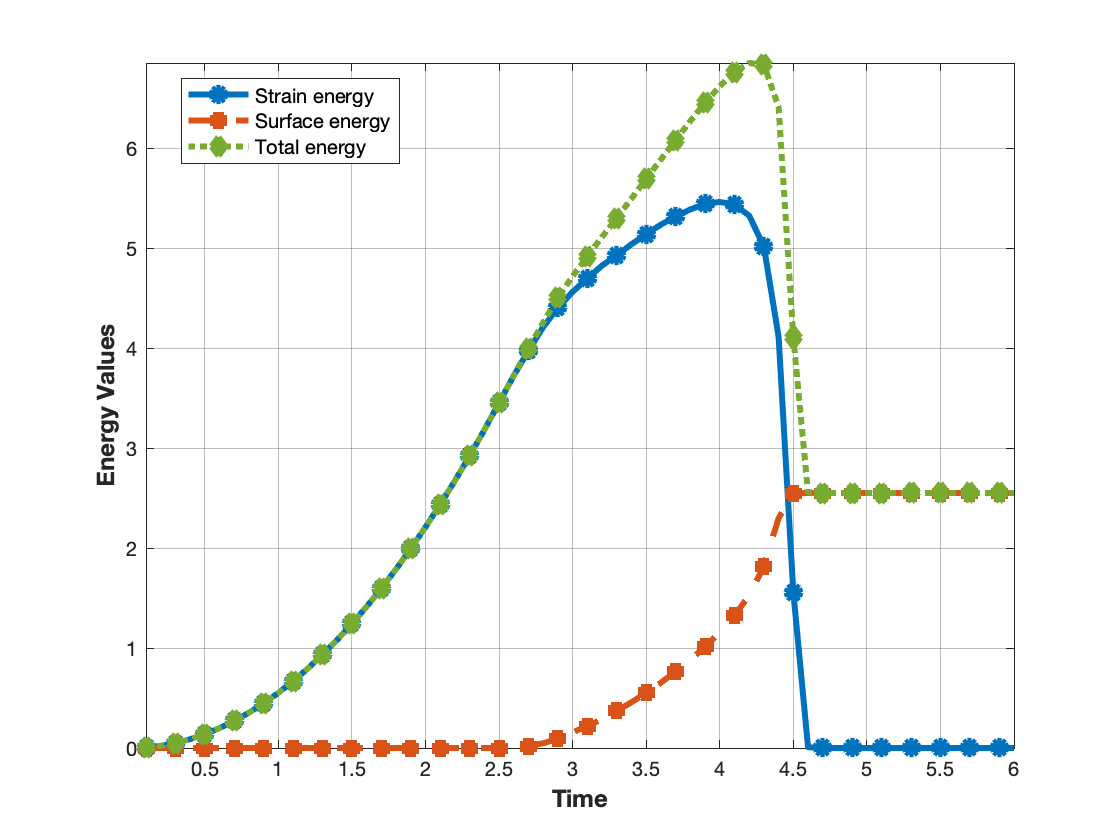}
\caption{Comparison of the energetic components (strain, surface, and total energy) for a model that uses a locally varying length scale parameter, $\xi$.}
    \label{fig:energies_LA_XI}
\end{figure}

Figure~\ref{fig:energies_LA_XI} shows the typical energy changes that happen when a material fractures in a phase-field simulation. It demonstrates how stored elastic energy turns into surface energy. At first, during the elastic loading from about $t=0$ to $t=0.25$, the system only stores energy through deformation. In this stage, the strain and total energy are the same and rise together, with no surface energy since no damage has occurred yet. Damage starts around $t=0.25$, indicated by an increase in surface energy. After this point, the energy put into the system is split between stored strain energy and the energy used to create new crack surfaces. This period of steady growth ends in a sudden failure when the strain energy hits its highest point at roughly $t=0.4$. This peak represents the material's maximum load capacity. These are very distinctive results compared with the classical linear elasticity damage model \cite{fernando2025xi}.  After the peak, the stored elastic energy quickly converts into surface energy, causing the strain energy to drop sharply and the crack to spread rapidly. This process leads to the final post-failure state for $t > 0.45$, where all the strain energy has been used up, and the system's total energy becomes constant, matching the final surface energy of the newly formed crack.


\section{Conclusion}\label{conclusions}

In this study, we introduced a novel phase-field framework for fracture simulation based on a modified Francfort-Marigo energy functional and its Ambrosio-Tortorelli (\textsf{AT1}-type) regularization. The proposed framework is within the context of an algebraically nonlinear, yet geometrically linear, strain-limiting theory of elasticity, an important distinction given its nonlinear strain energy density function. A key innovation of our approach lies in its three-field variational formulation, which concurrently solves for mechanics and phase-field variables, and a dynamically determined, spatially varying regularization parameter, $\xi(\bfx)$. The local fracture length scale is intrinsically governed by $\xi(\bfx)$. Its incorporation into the global energy minimization problem provides a robust theoretical basis for an adaptive meshing strategy that directly couples element size to the material's evolving physical behavior. We implemented this model using a custom solver built upon the open-source \textsf{deal.II} library and rigorously tested it on the challenging problem of quasi-static crack propagation under anti-plane shear.

Our numerical results demonstrate the efficacy of this formulation in driving localized mesh refinement. By judiciously reducing the length scale parameter in critical regions, our method significantly enhances both computational efficiency and approximation accuracy. This allows for the generation of high-fidelity solutions with substantially fewer degrees of freedom compared to conventional methods that necessitate globally fine meshes. The inherent capacity to concentrate computational resources positions our approach as a highly tractable and powerful tool. This work lays a strong foundation for future advancements, particularly in areas such as large-scale 3D simulations, dynamic fracture modeling, and other complex multi-physics applications where computational efficiency is paramount. 

\section*{Acknoledgement}
This material is based on work supported by the National Science Foundation under Grant No. 2316905.

\bibliographystyle{plain}
\bibliography{references}

\end{document}